\documentclass[12pt]{article}          
\usepackage{amssymb}
\usepackage{latexsym}
\usepackage{comment}
\usepackage{graphicx}
\usepackage{curves}
\usepackage{psfrag}  
\usepackage{color}

\newcommand{\K}{\mathcal{K}} 
\newcommand{\E}{\mathbb{E}^3}

\title{Polyhedra, Complexes, Nets and Symmetry}
\author{Egon Schulte\thanks{schulte@neu.edu}\\
Northeastern University\\
Department of Mathematics\\
Boston, MA 02115, USA}

\date{\sl\small\today }

\begin{document}

\maketitle

\begin{abstract}
Skeletal polyhedra and polygonal complexes in ordinary Euclidean 3-space are finite or infinite 3-periodic structures with interesting geometric, combinatorial, and algebraic properties. They can be viewed as finite or infinite 3-periodic graphs (nets) equipped with additional structure imposed by the faces, allowed to be skew, zig-zag, or helical. A  polyhedron or complex is {\em regular\/} if its geometric symmetry group is transitive on the flags (incident vertex-edge-face triples). There are 48 regular polyhedra (18 finite polyhedra and 30 infinite apeirohedra), as well as 25 regular polygonal complexes, all infinite, which are not polyhedra. Their edge graphs are nets well-known to crystallographers, and we identify them explicitly. There also are 6 infinite families of {\em chiral\/} apeirohedra, which have two orbits on the flags such that adjacent flags lie in different orbits. 
\end{abstract}

\section{Introduction}

Polyhedra and polyhedra-like structures in ordinary euclidean $3$-space $\mathbb{E}^3$ have been studied since the early days of geometry (Coxeter, 1973). However, with the passage of time, the underlying mathematical concepts and treatments have undergone fundamental changes. 

Over the past 100 years we can observe a shift from the classical approach viewing a polyhedron as a solid in space, to  topological approaches focussing on the underlying maps on surfaces (Coxeter \& Moser, 1980), to combinatorial approaches highlighting the basic incidence structure but deliberately suppressing the membranes that customarily span the faces to give a surface. 

These topological and combinatorial perspectives are already appearing in the well-known Kepler-Poinsot polyhedra and Petrie-Coxeter polyhedra (Coxeter, 1937, 1973). They underlie the skeletal approach to polyhedra proposed in (Gr\"unbaum, 1977a, 1977b) that has inspired a rich new theory of geometric polyhedra and symmetry (Dress, 1981, 1985; McMullen \& Schulte, 1997, 2002; McMullen, in prep.). 

The polygonal complexes described in this paper form an even broader class of discrete skeletal space structures than polyhedra. Like polyhedra they are comprised of vertices, joined by edges, assembled in careful fashion into polygons, the faces, allowed to be skew or infinite (Pellicer \& Schulte, 2010, 2013). However, unlike in polyhedra, more than two faces are permitted to meet at an edge. 

The regular polygonal complexes in $\mathbb{E}^3$ are finite structures or infinite $3$-periodic structures with crystallographic symmetry groups exhibiting interesting geometric, combinatorial and algebraic properties. This class includes the regular polyhedra but also many unfamiliar figures, once the planarity or finiteness of the polygonal faces is abandoned. Because of their skeletal structure, polygonal complexes are of natural interest to crystallographers. 

Regular polyhedra traditionally play a prominent role in crystal chemistry (Wells, 1977; O'Keeffe \& Hyde, 1996; Delgado-Friedrichs et al. 2005; O'Keeffe et al, 2008; O'Keeffe, 2008). There is considerable interest in the study of 3-periodic nets and their relationships to polyhedra. Nets are 3-periodic geometric graphs in space that represent crystal structures, in the simplest form with vertices corresponding to atoms and edges to bonds. The edge graphs of almost all regular polygonal complexes in $\E$ are highly-symmetric nets, with the only exceptions arising from the polyhedra which are not 3-periodic. We explicitly identify the nets by building on the work of (O'Keeffe, 2008) and exploiting the methods developed in (Delgado-Friedrichs et al., 2003).

Symmetry of discrete geometric structures is a frequently recurring theme in science. Polyhedral structures occur in nature and art in many contexts that a priori have little apparent relation to regularity (Fejes T\'oth, 1964; Senechal, 2013). Their occurrence in crystallography as crystal nets is a prominent example. See also (Wachman et al, 1974) for an interesting architecture inspired study of polyhedral structures that features a wealth of beautiful illustrations of figures reminiscent of skeletal polyhedral structures.

The present paper is organized as follows. In Section~\ref{pocone} we investigate basic properties of polygonal complexes, polyhedra, and nets, in particular focussing on structures with high symmetry. Section~\ref{symgroup} is devoted to the study of the symmetry groups of regular polygonal complexes as well as chiral polyhedra. In Sections~\ref{regpoly} and \ref{chipoly} we review the complete classification of the regular and chiral polyhedra following (McMullen \& Schulte, 2002, Ch. 7E) and (Schulte, 2004, 2005), respectively. The final Section~\ref{regpolycom} describes the classification of the regular polygonal complexes (Pellicer \& Schulte, 2010, 2013). Along the way we determine the nets of all 3-periodic regular polygonal complexes.

\section{Polyhedra, Complexes and Nets}
\label{pocone}

Polygonal complexes are geometric realizations in $\mathbb{E}^3$ of abstract incidence complexes of rank $3$ with polygon faces, that is, of abstract polygonal complexes (Danzer \& Schulte, 1982). We elaborate on this aspect in Section~\ref{polygcom}. Here we will not require familiarity with incidence complexes. However, occasionally it is useful to bear this perspective in mind. Polyhedra are the polygonal complexes with just two faces meeting at an edge.  

\subsection{Polygonal complexes}
\label{polygcom}

Following (Gr\"unbaum, 1977a), a {\em finite polygon\/}, or an {\em $n$-gon\/} (with $n\geq 3$), consists of a sequence $(v_1, v_2, \dots, v_n)$ of $n$ distinct points in $\mathbb{E}^3$, as well as of the line segments $(v_1, v_2), (v_2,v_3), \ldots, (v_{n-1},v_n)$ and $(v_n, v_1)$. Note that we are not making a topological disc spanned into the polygon part of the definition of a polygon. In particular, unless stated otherwise, the term ``convex polygon" refers only to the boundary edge path of what is usually called a convex polygon; that is, we ignore the interior. 

A (discrete) {\em infinite polygon\/}, or {\em apeirogon\/}, similarly consists of an infinite sequence of distinct points $(\dots, v_{-2},v_{-1}, v_0, v_1, v_2, \dots)$ in $\mathbb{E}^3$, as well as of the line segments $(v_i, v_{i+1})$ for each $i$, such that each compact subset of $\mathbb{E}^3$ meets only finitely many line segments. 

In either case the points are the {\em vertices\/} and the line segments the {\em edges\/} of the polygon.

Following (Pellicer \& Schulte, 2010), a {\em polygonal complex}, or simply a {\em complex}, $\K$ in $\mathbb{E}^{3}$ is a triple $(\mathcal{V},\mathcal{E},\mathcal{F}$) consisting of a set $\mathcal{V}$ of points, called {\em vertices}, a set $\mathcal{E}$ of line segments, called {\em edges}, and a set $\mathcal{F}$ of polygons, called {\em faces}, satisfying the following properties.
\begin{itemize}
\item[(a)] The graph $(\mathcal{V},\mathcal{E})$, the {\em edge graph\/} of $\K$, is connected.
\item[(b)] The vertex-figure of $\K$ at each vertex of $\K$ is connected. By the {\em vertex-figure\/} of $\K$ at a vertex $v$ we mean the graph, possibly with multiple edges, whose vertices are the vertices of $\K$ adjacent to $v$ and whose edges are the line segments $(u,w)$, where $(u, v)$ and $(v, w)$ are edges of a common face of $\K$. (There may be more than one such face in $\K$, in which case the corresponding edge $(u,w)$ of the vertex-figure at $v$ has multiplicity given by the number of such faces.)
\item[(c)] Each edge of $\K$ is contained in exactly $r$ faces of $\K$, for a fixed number $r \geq 2$.
\item[(d)] $\K$ is {\em discrete\/}, in the sense that each compact subset of $\mathbb{E}^{3}$ meets only finitely many faces of $\K$.
\end{itemize}

A (geometric) {\em polyhedron\/} in $\mathbb{E}^{3}$ is a polygonal complex with $r=2$. Thus, a polyhedron is a complex in which each edge lies in exactly two faces. The vertex-figures of a polyhedron are finite (simple) polygonal cycles. An infinite polyhedron in $\E$ is also called an {\em apeirohedron}.

\begin{figure}
\label{petcube}
\centering
\begin{picture}(110,100)
\put(0,1){\thicklines\color{red}\line(1,0){75}}
\put(76,0){\thicklines\color{red}\line(0,1){75}}
\put(75,76){\thicklines\color{red}\line(3,2){30}}
\put(105,96){\thicklines\color{red}\line(-1,0){75}}
\put(31,95){\thicklines\color{red}\line(0,-1){75}}
\put(30,21){\thicklines\color{red}\line(-3,-2){30}}
\put(0,-1){\thicklines\color{blue}\line(1,0){75}}
\put(75,1){\thicklines\color{blue}\line(3,2){30}}
\put(106,20){\thicklines\color{blue}\line(0,1){75}}
\put(105,94){\thicklines\color{blue}\line(-1,0){75}}
\put(30,94){\thicklines\color{blue}\line(-3,-2){30}}
\put(1,75){\thicklines\color{blue}\line(0,-1){75}}
\put(0,74){\thicklines\color{green}\line(1,0){75}}
\put(75,74){\thicklines\color{green}\line(3,2){30}}
\put(104,20){\thicklines\color{green}\line(0,1){75}}
\put(105,21){\thicklines\color{green}\line(-1,0){75}}
\put(30,19){\thicklines\color{green}\line(-3,-2){30}}
\put(-1,75){\thicklines\color{green}\line(0,-1){75}}
\put(74,0){\thicklines\color{black}\line(0,1){75}}
\put(0,76){\thicklines\color{black}\line(1,0){75}}
\put(30,96.5){\thicklines\color{black}\line(-3,-2){30}}
\put(29,95){\thicklines\color{black}\line(0,-1){75}}
\put(105,19){\thicklines\color{black}\line(-1,0){75}}
\put(75,-1.5){\thicklines\color{black}\line(3,2){30}}
\multiput(0,0)(75,0){2}{\circle*{4}}
\multiput(0,0)(0,75){2}{\circle*{4}}
\put(75,75){\circle*{4}}
\multiput(30,20)(75,0){2}{\circle*{4}}
\multiput(30,20)(0,75){2}{\circle*{4}}
\put(105,95){\circle*{4}}
\end{picture}
\caption{The Petrie dual of the cube. Shown are its four faces (in red, blue, green and black). The faces are the Petrie polygons of the cube.}
\end{figure}
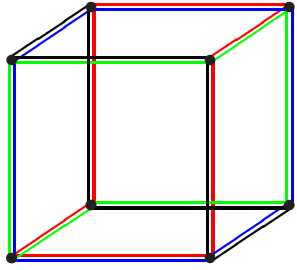

A simple example of a polyhedron is shown in Figure~\ref{petcube}. This is the ``Petrie-dual" of the cube obtained from the ordinary cube by replacing the square faces with the Petrie polygons while retaining the vertices and edges. Recall here that a {\em Petrie polygon\/}, or {\em $1$-zigzag\/}, of a polyhedron is a path along the edges such that any two, but no three, consecutive edges belong to a common face. The Petrie polygons of the cube are skew hexagons, and there are four of them. Hence the polyhedron in Figure~\ref{petcube} has 8 vertices, 12 edges, and 4 skew hexagonal faces, with 3 faces coming together at each vertex. 

Polyhedra are the best studied class of polygonal complexes and include the traditional convex polyhedra and star-polyhedra in $\mathbb{E}^{3}$ (Coxeter, 1973; Gr\"unbaum, 1977a; McMullen \& Schulte, 2002; McMullen, in prep.). When viewed purely combinatorially, the set of vertices, edges, and faces of a geometric polyhedron, ordered by inclusion, is an {\em abstract polyhedron\/}, or an {\em abstract polytope\/} of rank $3$. More generally, the underlying combinatorial ``complex" determined by the vertices, edges, and faces of any polygonal complex $\K$ (given by the triple $(\mathcal{V},\mathcal{E},\mathcal{F}$)), ordered by inclusion, is an {\em incidence complex\/} of rank $3$ in the sense of (Danzer \& Schulte, 1982). Here the term ``rank" refers to the ``combinatorial dimension" of the object; thus, rank $3$ complexes are incidence structures made up of 
objects called vertices (of rank~$0$), edges (of rank~$1$), and faces (of rank~$2$), in a purely combinatorial sense.

\begin{figure}
\label{skel2}
\centering
\begin{picture}(130,125)
\multiput(7.5,0)(0,30){4}{
\begin{picture}(120,30)
\thicklines
\multiput(0,0)(30,0){4}{\circle*{2}}
\multiput(12.5,8.33)(30,0){4}{\circle*{2}}
\multiput(25,16.66)(30,0){4}{\circle*{2}}
\multiput(0,0)(12.5,8.33){4}{\line(1,0){90}}
\multiput(0,0)(30,0){4}{\line(3,2){37.5}}
\end{picture}}
\put(7.5,0){
\begin{picture}(120,110)
\thicklines
\multiput(0,0)(30,0){4}{\line(0,1){90}}
\multiput(12.5,8.33)(30,0){4}{\line(0,1){90}}
\multiput(25,16.66)(30,0){4}{\line(0,1){90}}
\multiput(37.5,25)(30,0){4}{\line(0,1){90}}
\end{picture}}
\thicklines
\put(52,39.5){\color{red}{\line(3,2){12.5}}}
\put(52,38){\color{red}{\line(3,2){12.5}}}
\put(52,39.5){\color{red}{\line(3,2){12.5}}}
\put(22,38){\color{red}{\line(1,0){60}}}
\put(22,38){\color{red}{\line(3,2){12.5}}}
\put(82,38){\color{red}{\line(3,2){12.5}}}
\put(34.5,46.33){\color{red}{\line(1,0){60}}}
\put(52,8){\color{red}{\line(0,1){60}}}
\put(22,38){\color{red}{\line(3,2){12.5}}}
\put(52,68){\color{red}{\line(3,2){12.5}}}
\put(52,8){\color{red}{\line(3,2){12.5}}}
\put(64.5,16.33){\color{red}{\line(0,1){60}}}
\end{picture}
\caption{The 2-skeleton of the cubical tessellation. Each edge lies in four square faces.}
\end{figure}
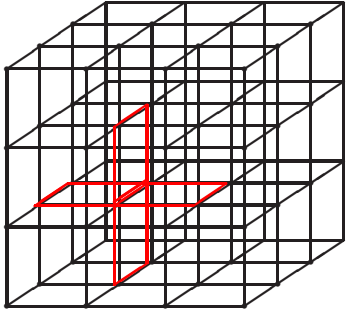

An easy example of an infinite polygonal complex in $\mathbb{E}^{3}$ which is not a polyhedron is given by the vertices, edges, and square faces of the standard cubical tessellation $\mathcal{C}$ (see Figure~\ref{skel2}). This complex $\K$ is called the {\em $2$-skeleton\/} of $\mathcal{C}$; each edge lies in four square faces so $r=4$. The tiles (cubes) of $\mathcal{C}$ are irrelevant in this context. A finite complex with $r=3$ can similarly be derived from the $2$-skeleton of the $4$-cube projected (as a Schlegel diagram) into $\mathbb{E}^{3}$. 

The cubical tessellation $\mathcal{C}$ in $\E$ gives rise to several other interesting complexes. For example, the family of all Petrie polygons of all cubes in $\mathcal{C}$ gives the (hexagonal) faces of a polygonal complex in which every edge belongs to exactly eight faces; the vertices and edges are just those of $\mathcal{C}$. This is the complex $\K_{6}(1,2)$ appearing in Table~\ref{tabsimply} below. Note here that every edge of a cube in $\mathcal{C}$ belongs to precisely two Petrie polygons of this cube; since every edge of $\mathcal{C}$ belongs to four cubes in $\mathcal{C}$, every edge of $\mathcal{C}$ must belong to eight Petrie polygons of cubes so $r=8$. 

If only the Petrie polygons of alternate cubes in $\mathcal{C}$ are taken as faces, we obtain a ``subcomplex" in which every edge belongs to exactly four hexagonal faces so then $r=4$.

The set of vertices, edges, and {\em triangular\/} faces of a cuboctahedron is not a polygonal complex; each edge lies in only one face, and the vertex-figures are not connected. Although our terminology could be adapted to cover polygonal structures in which some edges lie only in one face (as in this last example), we will explicitly exclude them here. Moreover, as we are mainly interested in highly symmetric structures, our definition includes the homogeneity condition (c). This condition is automatically satisfied for any polygonal structures with sufficiently high symmetry (for example, as given by the edge-transitivity of the symmetry group), provided at least two faces meet at an edge.

However, for an investigation of polygonal structures regardless of symmetry it is useful to replace part (c) in the definition of a polygonal complex by the following weaker requirement:

\smallskip
\begin{itemize}
\item[(c')] Each edge of $\K$ is contained in at least two faces of $\K$.
\end{itemize}
\smallskip

\noindent In this paper we will not require any of these modifications.

\subsection{Highly symmetric complexes}

There are several distinguished classes of highly symmetric polygonal complexes, each characterized by a distinguished transitivity property of the symmetry group. Some of these classes have analogues in the traditional theory of polyhedra but others feature characteristics that do not occur in the classical theory. 

We let $G(\K)$ denote the symmetry group of a polygonal complex $\K$, that is, the group of all Euclidean isometries of the affine hull of $\K$ that map $\K$ to itself. (Except when $\K$ is planar, this 
affine hull is $\mathbb{E}^{3}$ itself.) 

The most highly symmetric polygonal complexes $\K$ are those that we call {\em regular\/}, meaning that the symmetry group $G(\K)$ is transitive on the flags of $\K$. A {\em flag} of $\K$ is an incident triple consisting of a vertex, an edge and a face of $\K$. 

Two flags of $\K$ are called {\em $j$-adjacent\/} if they differ precisely in their elements of rank $j$, that is, their vertices, edges, or faces if $j=0$, $1$ or $2$, respectively. Flags are $j$-adjacent to only one flag if $j=0$ or $1$, or precisely $r-1$ flags if $j=2$. For example, in the $2$-skeleton of the cubical tessellation (with $r=4$) shown in Figure~\ref{skel2} every flag has exactly three $2$-adjacent flags. For polyhedra, every flag has one $j$-adjacent flag for every rank~$j$. 

The faces of a regular polygonal complex are (finite or infinite) congruent {\em regular\/} polygons in $\E$, with ``regular" meaning that their geometric symmetry group is transitive on the flags of the polygon. (A flag of a polygon is an incident vertex-edge pair.) 

Note that regular polygons in $\mathbb{E}^{3}$ are necessarily of one of the following kinds:\ finite, planar (convex or star-) polygons or non-planar ({\em skew\/}) polygons; (infinite) {\em apeirogons\/}, either planar zigzags or helical polygons; or linear polygons, either a line segment or a linear apeirogon with equal-sized edges (Gr\"unbaum, 1977a; Coxeter, 1991). We can show that linear regular polygons do not occur as faces of polygonal complexes.

We call a polygonal complex $\K$ {\em semiregular\/} (or {\em uniform\/}) if the faces of $\K$ are regular polygons (allowed to be non-planar or infinite) and $G(\K)$ is transitive on the vertices. 

A polygonal complex $\K$ is said to be {\em vertex-transitive\/}, {\em edge-transitive\/}, or {\em face-transitive} if $G(\K)$ is transitive on the vertices, edges, or faces, respectively. A complex which is vertex-transitive, edge-transitive, and face-transitive is called {\em totally transitive\/}. Every regular complex is totally transitive, but not vice versa. We call $\K$ a {\em $2$-orbit\/} polygonal complex if $\K$ has precisely two flag orbits under the symmetry group (Hubard, 2010; Cutler \& Schulte, 2011).

\subsection{Chiral polyhedra} 

Chiral polyhedra are arguably the most important class of $2$-orbit complexes. A (geometric) polyhedron $\K$ is {\em chiral\/} if $G(\K)$ has exactly two orbits on the flags such that any two adjacent flags are in distinct orbits (Schulte, 2004, 2005). This notion of chirality for polyhedra is different from the standard notion of chirality used in crystallography, but is inspired by it. The proper setting is that of a ``chiral realization" in $\mathbb{E}^{3}$ of an abstractly chiral or regular abstract polyhedron, where here abstract chirality or regularity are defined as above, but now in terms of the combinatorial automorphism group of the abstract polyhedron, not the geometric symmetry group. In a  sense that can be made precise, abstract chiral polyhedra occur in a ``left-handed" and a ``right-handed" version (Schulte \& Weiss, 1991, 1994), although this handedness is combinatorial and not geometric. 

Thus a chiral geometric polyhedron has maximum symmetry by ``combinatorial rotation" (but not by ``combinatorial reflection"), and has all its ``rotational"  combinatorial symmetries realized by euclidean isometries (but not in general by euclidean rotations). A regular geometric polyhedron has maximum symmetry by ``combinatorial reflection", and has all its combinatorial symmetries realized by euclidean isometries. 

\subsection{Nets}

Nets are important tools used in the modeling of $3$-periodic structures in crystal chemistry and materials science. A {\it net\/} in $\mathbb{E}^3$ is a $3$-periodic connected (simple) graph with straight edges (Delgado-Friedrichs \& O'Keeffe, 2005; Wells, 1977). Recall here that a figure in $\mathbb{E}^3$ is said to be {\em $3$-periodic\/} if the translation subgroup of its symmetry group is generated by translations in 3 independent directions. 

Highly symmetric nets often arise as the edge graphs (the graphs formed by the vertices and edges) of $3$-periodic higher rank structures in $\E$ such as $3$-dimensional tilings, apeirohedra, and infinite polygonal complexes. Note that the combinatorial automorphism group of a net $\mathcal{N}$ can be larger than its geometric symmetry group $G(\mathcal{N})$. The underlying abstract infinite graph of a net can often be realized in several different ways as a net in $\E$, and there is a natural interest in finding the maximum symmetry realization of this graph as a net. 

The Reticular Chemistry Structure Resource (RCSR) database located at the website http://rcsr.anu.edu.au contains a rich collection of crystal nets including in particular the most symmetric examples (O'Keeffe et al., 2008). It has become common practice to denote a net by a bold-faced 3-letter symbol such as \textbf{abc}. Examples are described below. The symbol is often a short-hand for a ``famous" compound represented by the net, or for the finer geometry of the net. A net can represent many compounds. Newly discovered nets tend to be named after the person(s) who discovered the net (using initials etc.). Many nets have ``alternative symbols" but we will not require them here. 

While the RCSR database contains information about 2000 named nets, TOPOS is a more recent research tool for the geometrical and topological analysis of crystal structures with a database of more than 70,000 nets (Blatov, 2012); it is currently under further development by Vladislav Blatov and Davide Proserpio.

The nets $\mathcal{N}$ occurring in this paper are {\it uninodal\/}, meaning that $G(\mathcal{N})$ is transitive on the vertices (nodes) of $\mathcal{N}$. For a vertex $v$ of a net $\mathcal{N}$, the convex hull of the neighbors of $v$ in $\mathcal{N}$ is called the {\em coordination figure\/} of $\mathcal{N}$ at $v$. 

The edge graph of each regular polygonal complex $\K$ is a net referred to as the {\it net of the complex\/}. The identification of the nets arising as edge graphs of regular polygonal complexes is greatly aided by the fact that there already exists a classification of the nets in $\E$ that are called regular or quasiregular in the chemistry literature (Delgado-Friedrichs et al., 2003, 2005). Although this terminology for nets is not consistent with our terminology for polyhedral complexes, we will maintain it for the convenience of the reader. Note that a net of a regular complex may have symmetries which are not symmetries of the complex.

A net $\mathcal{N}$ in $\E$ is called {\it regular\/} if $\mathcal{N}$ is uninodal and if, for each vertex $v$ of $\mathcal{N}$, the coordination figure of $\mathcal{N}$ at $v$ is a regular convex polygon (in the ordinary sense) or a Platonic solid whose own rotation symmetry group is a subgroup of the stabilizer of $v$ in $G(\mathcal{N})$ (the {\em site symmetry group\/} of $v$ in $\mathcal{N}$). Here the rotation symmetry group of a regular convex polygon is taken relative to~$\E$ and is generated by two half-turns in~$\E$.

As pentagonal symmetry is impossible, the coordination figures of a regular net must necessarily be triangles, squares, tetrahedra, octahedra or cubes. 

A net $\mathcal{N}$ is called {\it quasiregular\/} if $\mathcal{N}$ is uninodal and the coordination figure of $\mathcal{N}$ at every vertex is a quasiregular convex polyhedron. Recall that a {\em quasiregular\/} convex polyhedron in $\E$ is a semiregular convex polyhedron (with regular faces and a vertex-transitive symmetry group) with exactly two kinds of faces alternating around each vertex. There are only two quasiregular convex polyhedra in $\E$ (Coxeter, 1973), namely the well-known cuboctahedron $3.4.3.4$ and icosidodecahedron $3.5.3.5$, of which the latter cannot occur because of its icosahedral symmetry. Thus the coordination figures of a quasiregular net are cuboctahedra. Conversely, a uninodal net with cuboctahedra as coordination figures is necessarily quasiregular.  

There are exactly five regular nets in $\E$, one per possible coordination figure. Following (O'Keeffe et al., 2008; Delgado-Friedrichs et al., 2003) these nets are denoted by {\bf srs}, {\bf nbo}, {\bf dia}, {\bf pcu} and {\bf bcu}; their coordination figures are triangles, squares, tetrahedra, octahedra or cubes, respectively. It is also known that there is just one quasiregular net in $\E$, denoted {\bf fcu}, with coordination figure a cuboctahedron. These nets have appeared in many publications, often under different names (O'Keeffe \& Hyde, 1996; Delgado Friedrichs et al, 2003).

The net {\bf srs} is the net of Strontium Silicide $\rm{SrSi}_2$, hence the notation. It coincides with the ``net (10,3)-a" of (Wells, 1977), the ``Laves net" of (Pearce, 1978), the ``Y$^*$ lattice complex" of (Koch \& Fischer, 1983), the net ``3/10/c1" of (Koch \& Fischer, 1995), as well as the ``labyrinth graph of the gyroid surface" (Hyde \& Ramsden, 2000). For an account on the history of the {\bf srs} net and its gyroid surface see also (Hyde et al., 2008). The notation {\bf nbo} signifies the net of Niobium Monoxide NbO and coincides with the lattice complex J$^*$ of (Koch \& Fischer, 1983). The net {\bf dia} is the famous diamond net, or ``lattice complex D" (Koch \& Fischer, 1983), which is the net of the diamond form of carbon. The names {\bf pcu}, {\bf fcu}, and {\bf bcu} stand for the ``primitive cubic lattice" (the standard cubic lattice), the ``face-centered cubic lattice", and the ``body-centered" cubic lattice in $\E$, respectively; these are also known as the lattice complexes $c$P, $cI$ and $c$F, respectively. 

In addition we will also meet the nets denoted {\bf acs}, {\bf sod}, {\bf crs} and {\bf hxg} (O'Keeffe et al., 2008; O'Keeffe, 2008). The net {\bf acs} is named after Andrea C. Sudik (Sudik et al., 2005) and is observed in at least 177 compounds (according to the TOPOS database). The net {\bf sod\/} represents the sodalite structure and can be viewed as the edge-graph of the familiar tiling of $\E$ by truncated octahedra also known as the Kelvin-structure (Delgado-Friedrichs et al., 2005). The symbol {\bf crs} labels the net of the oxygens coordination in idealized beta-cristobalite, but also represents at least 10 other compounds (according to the TOPOS database); the net is also known as {\bf dia-e} and 3d-kagom\`{e}, and appears as the 6-coordinated sphere packing net corresponding to the cubic invariant lattice complex T (Koch \& Fischer, 1983). Finally, {\bf hxg\/} has a regular hexagon as its coordination figure; its ``augmented" net (obtained by replacing the original vertices by the edge graph of the coordination figures) gives a structure representing polybenzene. 

\section{Symmetry Groups of Complexes and Polyhedra}
\label{symgroup}

The symmetry groups $G:=G(\K)$ hold the key to the structure of regular polygonal complexes $\K$. They have a distinguished generating set of subgroups $G_0,G_1,G_2$ obtained as follows. Via a variant of Wythoff's construction (Coxeter, 1973), these subgroups enable us to recover a regular complex from its group. 

\subsection{Distinguished generators}

Choose a fixed, or {\em base}, flag $\Phi := \{F_0, F_1, F_2\}$ of $\K$, where $F_{0},F_{1},F_{2}$, respectively, denote the vertex, edge, or face of the flag. For $i=0,1,2$ let $G_i$ denote the stabilizer of $\Phi \setminus \{F_i\}$ in $G$; this is the subgroup of $G$ stabilizing every element of $\Phi$ except $F_i$. For example, $G_2$ consists of all symmetries of $\K$ fixing $F_0$ and $F_1$ and thus fixing the entire line through $F_1$ pointwise.  

Also, for $\Psi\subseteq\Phi$ define $G_{\Psi}$ to be the stabilizer of $\Psi$ in $G$, that is, the subgroup of $G$ stabilizing every element of $\Psi$. Then $G_i$ is just $G_{\{F_j, F_k\}}$ for each $i$, where here $\{i, j, k\}=\{0,1,2\}$. Moreover, $G_\Phi$ is the stabilizer of the base flag $\Phi$ itself. We also write $G_{F_i}:=G_{\{F_i\}}$ for $i=0,1,2$; this is the stabilizer of $F_i$ in $G$. 

The subgroups $G_0,G_1,G_2$ have remarkable properties. In particular,
\[ G = \langle G_0,G_1,G_2 \rangle \]
and
\[ G_0 \cap G_1 = G_0 \cap G_2 = G_1 \cap G_2 = G_0 \cap G_1 \cap G_2 = G_\Phi .\]

While these properties already hold at the abstract level of incidence complexes (Schulte, 1983), the euclidean geometry of 3-space comes into play when we investigate the possible size of the flag stabilizer $G_\Phi$. It turns out that there are two possible scenarios. 

In fact, either $G_{\Phi}$ is trivial, and then the (full) symmetry group $G$ is simply flag-transitive; or $G_{\Phi}$ has order $2$, the complex $\K$ has planar faces, and $G_{\Phi}$ is generated by the reflection in the plane of $\mathbb{E}^{3}$ containing the base face $F_2$ of $\K$. In the former case we call $\K$ a {\em simply flag-transitive} complex. In the latter case we say that $\K$ is {\em non-simply flag-transitive}, or that $\K$ has {\em face mirrors\/} since then the planes in $\E$ through faces of $\K$ are mirrors (fixed point sets) of reflection symmetries. 

If $\K$ is a simply flag-transitive complex, then $G_0 = \langle R_0 \rangle$ and $G_1 = \langle R_1 \rangle$, for some point, line or plane reflection $R_0$ and some line or plane reflection $R_1$; and $G_2$ is a cyclic or dihedral group of order $r$. (A reflection in a line is a half-turn about the line.) Moreover, 
\[ G_{F_0}=\langle R_1, G_2 \rangle,\; G_{F_1}=\langle R_0, G_2 \rangle, \;
G_{F_2}=\langle R_0, R_1 \rangle\cong D_p ,\]
where $p$ is the number of vertices in a face of $\K$ and $D_p$ denotes the dihedral group of order $2p$ (allowing $p=\infty$). Since $\K$ is discrete, the stabilizer $G_{F_0}$ of $F_0$ in $G$ is necessarily a finite group called the {\em vertex-figure group\/} of $\K$ at $F_0$. In particular, the vertex-figure group $G_{F_0}$ acts simply flag-transitively on the finite graph forming the vertex-figure of $\K$ at $F_0$. (A flag of a graph is just an incident vertex-edge pair.) 

If $\K$ is a non-simply flag-transitive complex, then 
\[ G_0 \cong C_2 \times C_2 \cong G_1, \; G_2\cong D_r \]
and the vertex-figure group $G_{F_0}=\langle G_1,G_2 \rangle$ is again finite as $\K$ is discrete.

\subsection{The case of polyhedra}

The theory is particularly elegant in the case when $\K$ is a regular polyhedron. Then $\K$ is necessarily simply flag-transitive and $G_2$ is also generated by a reflection $R_2$ in a point, line or plane. Thus 
\[ G=\langle R_0,R_1,R_2\rangle \] 
and $G$ is a discrete (generalized) reflection group in $\mathbb{E}^{3}$, where here the term ``reflection group" refers to a group generated by reflections in points, lines or planes. 

A regular polyhedron has a (basic) Schl\"afli type $\{p,q\}$, where as above $p$ is the number of vertices in a given face and $q$ denotes the number of faces containing a given vertex; here $p$, but not $q$, may be $\infty$. The distinguished involutory generators $R_0,R_1,R_2$ of $G$ satisfy the standard Coxeter-type relations
\[ R_{0}^{2} = R_{1}^{2} = R_{2}^{2} = (R_{0}R_{1})^{p} = (R_{1}R_{2})^{q} = (R_{0}R_{2})^{2} = 1, \]
but in general there are other independent relations too; these additional relations are determined by the cycle structure of the edge graph, or net, of $\K$ (McMullen \& Schulte, 2002, Ch.~7E). 

For a chiral polyhedron $\K$, the symmetry group $G$ has two non-involutory {\em distinguished generators\/} $S_1,S_2$ obtained as follows. Let again $\Phi := \{F_{0},F_{1},F_{2}\}$ be a base flag, let $F_{0}'$ be the vertex of $F_1$ distinct from $F_0$, let $F_{1}'$ be the edge of $F_2$ with vertex $F_0$ distinct from $F_1$, and let $F_{2}'$ be the face containing $F_1$ distinct from $F_2$. Then the generator $S_{1}$ stabilizes the base face $F_2$ and cyclically permutes the vertices of $F_2$ in such a manner that $F_{1}S_{1} = F'_{1}$ (and $F'_{0}S_{1} = F_{0}$), while $S_{2}$ fixes the base vertex $F_0$ and cyclically permutes the vertices in the vertex-figure at $F_0$ in such a way that $F_{2}S_{2} = F'_{2}$ (and $F'_{1}S_{2} = F_{1}$). These generators $S_1,S_2$ satisfy (among others) the relations
\begin{equation}
\label{reltwo}
S_{1}^p = S_{2}^q = (S_{1}S_{2})^{2}  = 1, 
\end{equation}
where again $\{p,q\}$ is the Schl\"afli type of $\K$. 

Two alternative sets of generators for $G$ are given by $\{S_1,T\}$ and $\{T,S_2\}$, where $T:= S_{1}S_{2}$ is the involutory symmetry of $\K$ that interchanges simultaneously the two vertices of the base edge $F_1$ and the two faces meeting at $F_1$. Combinatorially speaking, $T$ acts like a half-turn about the midpoint of the base edge (but geometrically, $T$ may not be a half-turn about a line).

\subsection{Wythoff's construction}

Regular polygonal complexes and chiral polyhedra can be recovered from their symmetry groups $G$ by variants of the classical {\em Wythoff construction\/} (Coxeter, 1973). Two variants are needed, one essentially based on the generating subgroups $G_0,G_1,G_2$ of $G$ and applying only to polygonal complexes which are regular (McMullen \& Schulte, 2002, Ch. 5; Pellicer \& Schulte, 2010), and the other based on the generators $S_1,S_2$ and applying to both regular and chiral polyhedra. (To subsume regular polyhedra under the latter case it is convenient to set $S_{1}:=R_{0}R_{1}$, $S_{2}:=R_{1}R_{2}$ and $T:=S_{1}S_{2}=R_{0}R_{2}$.)

For regular polygonal complexes (including polyhedra), Wythoff's construction proceeds from the base vertex $v:=F_0$, the {\em initial\/} vertex, and builds the complex (or polyhedron) $\K$ as an orbit structure, beginning with the construction of the base flag. Relative to the set of generating subgroups, the essential property of $v$ is that it is invariant under all generating subgroups but the first; that is, $v$ is invariant under $G_{1}$ and $G_{2}$ but not $G_0$. The base vertex, $v$, is already given. The vertex sets of the base edge and base face of $\K$ are the orbits of $v$ under the subgroups $G_0$ and $\langle G_{0},G_{1}\rangle$, respectively. This determines the base edge as the line segment joining its vertices, and the base face as an edge path joining its vertices in succession. Once the base flag has been constructed, we simply obtain the vertices, edges and faces of $\K$ as the images under $G$ of the base vertex, base edge or base face, respectively.

For chiral polyhedra $\K$ we can similarly proceed from the alternative generators $T,S_2$ of $G$, again choosing $v:=F_0$ as the initial (or base) vertex. Now the base edge and base face of $\K$ are given by the orbits of $v$ under $\langle T\rangle$ and $\langle S_1\rangle$, respectively; and as before, the vertices, edges and faces of $\K$ are just the images under $G$ of the base vertex, base edge or base face, respectively. 

In practice, Wythoff's construction is often applied to groups that ``look like" symmetry groups of regular complexes or chiral polyhedra. In fact, this approach then often enables us to establish the existence of such structures. A necessary condition in this case is the existence of a common fixed point of $G_1$ and $G_2$, which then becomes the initial vertex. 

\section{Regular Polyhedra}
\label{regpoly}

Loosely speaking there are 48 regular polyhedra in $\E$, up to similarity (that is, congruence and scaling). They comprise 18 finite polyhedra and 30 (infinite) apeirohedra. We follow the classification scheme described in (McMullen \& Schulte, 2002, Ch. 7E) and designate these polyhedra by generalized Schl\"afli symbols that usually are obtained by padding the basic symbol $\{p,q\}$ with additional symbols signifying specific information (such as extra defining relations for the symmetry group in terms of the distinguished generators).

\subsection{Finite polyhedra}
The finite regular polyhedra are all derived from the five Platonic solids:\ the tetrahedron $\{3,3\}$, the octahedron $\{3,4\}$, the cube $\{4,3\}$, the icosahedron $\{3,5\}$, and the dodecahedron $\{5,3\}$. In addition to the Platonic solids, there are the four regular star-polyhedra, also known as {\em Kepler-Poinsot polyhedra\/}: the great icosahedron $\{3,\frac{5}{2}\}$, the great stellated dodecahedron $\{\frac{5}{2},3\}$, the great dodecahedron $\{5,\frac{5}{2}\}$, and the small stellated dodecahedron$\{\frac{5}{2},5\}$. (The fractional entries $\frac{5}{2}$ indicate that the corresponding faces or vertex-figures are star-pentagons.) These nine examples are the {\em classical\/} regular polyhedra. 

The remaining nine finite regular polyhedra are the Petrie duals of the nine classical regular polyhedra. The {\em Petrie dual\/} of a regular polyhedron is (usually) a new regular polyhedron with the same vertices and edges, obtained by replacing the faces by the Petrie polygons. The Petrie dual of the Petrie dual of a regular polyhedron is the original polyhedron. For example, the four (skew hexagonal) Petrie polygons of the cube form the faces of the Petrie dual of the cube, which is usually denoted $\{6,3\}_4$ and is shown in Figure~\ref{petcube}.  (The suffix indicates the length of the Petrie polygon, $4$ in this case.) The underlying abstract polyhedron corresponds to a map with four hexagonal faces on the torus.

\subsection{Apeirohedra}

The 30 apeirohedra fall into three families comprised of the 6 planar, the 12 ``reducible", and the 12 ``irreducible" examples. Their symmetry groups are crystallographic groups. The use of the terms ``reducible" and ``irreducible" for apeirohedra is consistent with what we observe at the group level:\ the symmetry group is affinely reducible or affinely irreducible, respectively. In saying that a group of isometries of $\mathbb{E}^3$ is {\em affinely reducible\/}, we mean that there is a line $l$ in $\E$ such that the group permutes the lines parallel $l$; the group then also permutes the planes perpendicular to $l$.

The six planar apeirohedra can be disposed off quickly:\ they are just the regular tessellations 
\[ \{3,6\},\; \{6,3\},\; \{4,4\}\] 
in the plane $\mathbb{E}^2$ by triangles, hexagons, and squares, respectively, as well as their Petrie duals, 
\[ \{\infty,6\}_{3},\; \{\infty,3\}_{6},\; \{\infty,4\}_{4},\] 
which have zig-zag faces. 

\subsection{Blended apeirohedra}

The ``reducible" apeirohedra are {\em blends\/}, in the sense that they are obtained by ``blending" a plane apeirohedron with a linear polygon (that is, a line segment $\{\;\}$ or an apeirogon $\{\infty\}$) contained in a line perpendicular to the plane. Thus there are $6\times 2 =12$ such blends. 

The two projections of a blended apeirohedron onto its component subspaces recover the two original components, that is, the original plane apeirohedron as well as the line segment or apeirogon. 

\begin{figure}
\centering
\begin{picture}(240,80)
\multiput(0,0)(0,36){2}{
\begin{picture}(240,50)
\thinlines
\multiput(0,0)(40,0){5}{\circle*{2}}
\multiput(15,10)(40,0){5}{\circle*{2}}
\multiput(30,20)(40,0){5}{\circle*{2}}
\multiput(45,30)(40,0){5}{\circle*{2}}
\multiput(60,40)(40,0){5}{\circle*{2}}
%\multiput(75,50)(40,0){5}{\circle*{2}}
%\multiput(90,60)(40,0){5}{\circle*{2}}
\multiput(0,0)(15,10){5}{\line(1,0){160}}
\multiput(0,0)(40,0){5}{\line(3,2){60}}
\end{picture}}
\multiput(192,20)(0.3,0.92){50}{\color{blue}\circle*{1}}
\multiput(0,0)(30,20){2}{
\multiput(0,0)(80,00){2}{
\multiput(2,0)(0.3,0.92){50}{\color{blue}\circle*{1}}
\multiput(2,0)(0.8,0.72){50}{\color{blue}\circle*{1}}
\multiput(17,46)(0.8,-0.72){50}{\color{blue}\circle*{1}} 
\multiput(42,36)(0.3,-0.52){50}{\color{blue}\circle*{1}} 
\put(2,0){\color{blue}\circle*{4}}
\put(17,46){\color{blue}\circle*{4}}
\put(42,36){\color{blue}\circle*{4}}
\put(57,10){\color{blue}\circle*{4}}}}
\multiput(40,0)(80,0){2}{
\multiput(2,36)(0.8,-0.72){50}{\color{blue}\circle*{1}} 
\multiput(2,36)(0.3,-0.52){50}{\color{blue}\circle*{1}} 
\multiput(42,0)(0.3,0.92){50}{\color{blue}\circle*{1}} 
\multiput(17,10)(0.8,0.72){50}{\color{blue}\circle*{1}} 
\put(2,36){\color{blue}\circle*{4}}
\put(57,46){\color{blue}\circle*{4}}
\put(42,0){\color{blue}\circle*{4}}
\put(17,10){\color{blue}\circle*{4}}}         %}
\multiput(15,10)(30,20){2}{
\multiput(40,0)(80,00){2}{
\multiput(2,0)(0.3,0.92){50}{\color{blue}\circle*{1}}
\multiput(2,0)(0.8,0.72){50}{\color{blue}\circle*{1}}
\multiput(17,46)(0.8,-0.72){50}{\color{blue}\circle*{1}} 
\multiput(42,36)(0.3,-0.52){50}{\color{blue}\circle*{1}} 
\put(2,0){\color{blue}\circle*{4}}
\put(17,46){\color{blue}\circle*{4}}
\put(42,36){\color{blue}\circle*{4}}
\put(57,10){\color{blue}\circle*{4}}}
\multiput(0,0)(80,0){2}{
\multiput(2,36)(0.8,-0.72){50}{\color{blue}\circle*{1}} 
\multiput(2,36)(0.3,-0.52){50}{\color{blue}\circle*{1}} 
\multiput(42,0)(0.3,0.92){50}{\color{blue}\circle*{1}} 
\multiput(17,10)(0.8,0.72){50}{\color{blue}\circle*{1}} 
\put(2,36){\color{blue}\circle*{4}}
\put(57,46){\color{blue}\circle*{4}}
\put(42,0){\color{blue}\circle*{4}}
\put(17,10){\color{blue}\circle*{4}}}}
\end{picture}
\caption{The blend of the square tessellation with the line segment. The vertices lie in two parallel planes, and over each square of the original square tessellation lies one skew square (tetragon) of the blend.}
\label{blendline}
\end{figure}
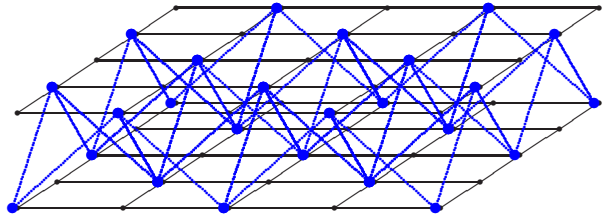

For example, the blend of the square tessellation in $\mathbb{E}^2$ with a line segment $[-1,1]$ positioned along the $z$-axis in $\E$ has its vertices in the planes $z=-1$ and $z=1$ parallel to $\mathbb{E}^2$, and is obtained from the square tessellation $\{4,4\}$ in $\mathbb{E}^2$ by alternately raising or lowering (alternate) vertices (see Figure~\ref{blendline}). Its faces are tetragons (skew squares), with vertices alternating between the two planes. Its designation is $\{4,4\}\#\{\;\}$, where $\#$ indicates the blending operation. 

The blend of the square tessellation in $\mathbb{E}^2$ with a linear apeirogon positioned along the $z$-axis is more complicated. Its faces are helical polygons rising in two-sided infinite vertical towers above the squares of the tessellation in such a way that the helical polygons over adjacent squares have opposite orientations (left-handed or right-handed) and meet along every fourth edge as they spiral around the towers. The designation in this case is $\{4,4\}\#\{\infty\}$. 

Strictly speaking, each blended apeirohedron belongs to an infinite family of apeirohedra obtained by (relative) rescaling of the two components of the blend; that is, each blended regular apeirohedron really represents a one-parameter family (of mutually non-similar) regular apeirohedra with the same combinatorial characteristics.

\begin{figure}
\label{blend44inf}
\centering
\begin{picture}(120,135)
\put(0,10){
\multiput(7.5,-10)(0,30){5}{
\begin{picture}(120,130)
\thinlines
\multiput(0,0)(30,0){3}{\circle*{1}}
\multiput(12.5,8.33)(30,0){3}{\circle*{1}}
\multiput(25,16.66)(30,0){3}{\circle*{1}}
\multiput(0,0)(12.5,8.33){3}{\line(1,0){60}}
\multiput(0,0)(30,0){3}{\line(3,2){25}}
\end{picture}}
\put(7.5,-10){
\begin{picture}(120,130)
\thinlines
\multiput(0,0)(30,0){3}{\line(0,1){120}}
\multiput(12.5,8.33)(30,0){3}{\line(0,1){120}}
\multiput(25,16.66)(30,0){3}{\line(0,1){120}}
\end{picture}}
\put(7.5,-10){
\begin{picture}(100,100)
\thicklines
\put(30,30){\color{red}{\line(1,1){30}}}
\put(60,60){\color{red}{\line(1,3){12.5}}}
\put(72.5,98.33){\color{red}{\line(-1,1){30}}}
\multiput(30,31)(0.25,-0.423){51}{\color{red}{\circle*{.007}}}
\put(30,30){\color{blue}{\line(-1,1){30}}}
\multiput(30,29)(0.25,-0.423){50}{\color{blue}{\circle*{.007}}}
\put(0,60){\color{blue}{\line(1,3){12.5}}}
\put(12.5,99){\color{blue}{\line(1,1){30}}}
\put(72.5,97.33){\color{red}{\line(-1,1){30}}}
\put(72.5,99.33){\color{green}{\line(-1,1){30}}}
\multiput(72.5,98.33)(0.25,-0.423){50}{\color{green}{\circle*{.007}}}
\put(85,76.66){\color{green}{\line(-1,-1){30}}}
\put(42.5,7){\color{green}{\line(1,3){12.5}}}
\put(42.5,9.66){\color{black}{\line(1,3){12.5}}}
\put(55,46.66){\color{black}{\line(-1,1){30}}}
\multiput(25,76.66)(-0.25,0.423){50}{\color{black}{\circle*{.007}}}
\put(12.5,97.2){\color{black}{\line(1,1){30}}}
\put(42.5,8.33){\circle*{3}}
\put(30,30){\circle*{3}}
\put(55,46.33){\circle*{3}}
\put(0,60){\circle*{3}}
\put(60,60){\circle*{3}}
\put(25,76.33){\circle*{3}}
\put(85,76.33){\circle*{3}}
\put(12.5,98.33){\circle*{3}}
\put(72.5,98.33){\circle*{3}}
\put(42.5,128.33){\circle*{3}}
\put(14,2){\color{blue}\scriptsize{\rm{blue}}}
\put(44,2){\color{red}\scriptsize{\rm{red}}}
\put(24,10.33){\color{black}\scriptsize{\rm{black}}}
\put(54,10.33){\color{green}\scriptsize{\rm{green}}}
\end{picture}} }
\end{picture}
\caption{The four helical facets of the blended apeirohedron $\{4,4\}\#\{\infty\}$ that share a vertex. Each vertical column over a square of the square tessellation is occupied by exactly one facet of $\{4,4\}\#\{\infty\}$, spiraling around the column.}
\end{figure}
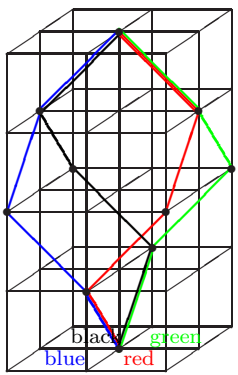

The six blends of a planar apeirohedron with a line segment $\{\;\}$ are listed in Table~\ref{blendlineseg}. They have their vertices in two parallel planes and hence are not $3$-periodic. Thus their edge graph cannot be a net. 

The edge graphs of the apeirohedra in the first two rows of Table~\ref{blendlineseg} are planar graphs, but they 
do not lie in a plane. In fact, the graphs are isomorphic to the edge graphs of the regular square tessellation $\{4,4\}$ or triangle tessellation $\{3,6\}$, respectively; this can be seen be projecting them onto the plane of the planar component of the blend. It follows that, under certain conditions, these apeirohedra have an edge graph that occurs as the contact graph of a sphere packing arrangement with spheres centered at the vertices. In fact, in order to have a faithful representation of the edge graph as a contact graph, the edges must be short enough to avoid forbidden contacts or overlaps between two spheres centered at vertices which are not joined by an edge. This condition depends on the relative scaling of the components of the blend but is satisfied if the line segment in the blend is short compared with the edge length of the planar apeirohedron in the blend (that is, if the two parallel planes containing the vertices of the blend are close to each other). The resulting sphere packings then are the sphere packings $4^{4}\rm{IV}$ and $3^{6}\rm{III}$ of (Koch \& Fischer, 1978, pp. 131-133).

On the other hand, the edge graphs of the apeirohedra in the third row of Table~\ref{blendlineseg} cannot be contact graphs of sphere packing arrangements (for any relative scaling of the components of the blend), as there always are forbidden overlaps.
\smallskip

\begin{table}
\label{blendlineseg} 
%\begin{center}
\centering
\begin{tabular}{lll}
Blend&  & Sphere Packing\\
\hline
$\{4,4\}\#\{\;\}$ & $\{\infty,4\}_{4}\#\{\;\}$ & $4^{4}\rm{IV}$ \\
$\{3,6\}\#\{\;\}$ & $\{\infty,6\}_{3}\#\{\;\}$ & $3^{6}\rm{III}$\\
$\{6,3\}\#\{\;\}$ & $\{\infty,3\}_{6}\#\{\;\}$ & \\
\end{tabular}
%\end{center}
\caption{The edge graphs of the blends of a planar regular apeirohedron with a line segment. Petrie dual apeirohedra have the same edge graph and are listed in the same row. Where applicable, the third column lists the corresponding sphere packing of (Koch \& Fischer,~1978).}
\end{table}

On the other hand, the blends with the linear apeirogon $\{\infty\}$ have $3$-periodic edge graphs and hence yield highly symmetric nets. These nets were already identified in (O'Keeffe, 2008) and are listed in Table~\ref{blendnet} using the notation for nets described earlier. 
\smallskip

\begin{table}
\label{blendnet} 
%\begin{center}
\centering
\begin{tabular}{lll}
Blend&  &Net \\
\hline
$\{4,4\}\#\{\infty\}$ & $\{\infty,4\}_{4}\#\{\infty\}$ & {\bf dia} \\
$\{3,6\}\#\{\infty\}$ & $\{\infty,6\}_{3}\#\{\infty\}$ & {\bf pcu}\\
$\{6,3\}\#\{\infty\}$ & $\{\infty,3\}_{6}\#\{\infty\}$ & {\bf acs}\\
\end{tabular}
%\end{center}
\caption{The nets of the blends of a planar regular apeirohedron with a linear apeirogon. Petrie dual apeirohedra have the same net and are listed in the same row.}
\end{table}

\subsection{Pure apeirohedra}

\begin{figure}
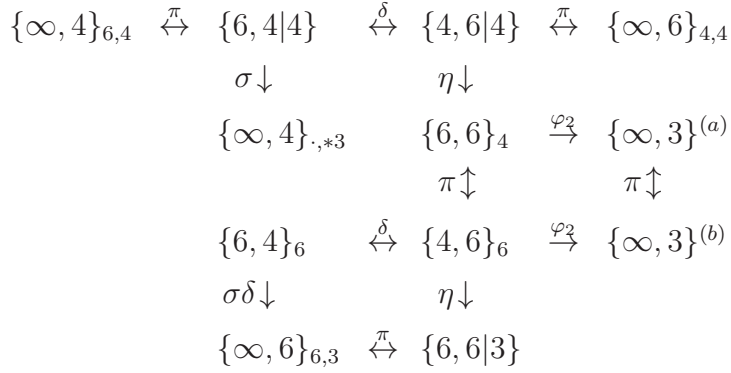

\centering
\label{diag}
$\begin{array}{lllllll}
\{\infty,4\}_{6,4} & \stackrel{\pi}{\leftrightarrow}\!\, & 
\{6,4|4\} & \!\stackrel{\delta}{\leftrightarrow}\! & \{4,6|4\} &
\stackrel{\pi}{\leftrightarrow} & \{\infty,6\}_{4,4}  \\[.07in]
& &\;\; \sigma\! \downarrow & &\;\; \eta\!\downarrow  & &   \\[.07in]
& & \{\infty,4\}_{\cdot,*3} & & \{6,6\}_{4} &
\stackrel{\varphi_2}{\rightarrow} & \{\infty,3\}^{(a)}  \\[.05in]
& & & &\;\; \pi\! \updownarrow & &\;\; \pi\!\updownarrow   \\[.07in] 
& & \{6,4\}_{6} & \!\stackrel{\delta}{\leftrightarrow}\! & \{4,6\}_{6}
& \stackrel{\varphi_2}{\rightarrow} & \{\infty,3\}^{(b)}  \\[.07in]
& & \,\sigma\delta\!\!\,\downarrow & & \;\;\eta\!\downarrow  & &   \\[.07in]
& & \{\infty,6\}_{6,3} & \!\stackrel{\pi}{\leftrightarrow}\! &
\{6,6 |3\} & & 
\end{array}$
\caption{Relationships among the twelve pure regular apeirohedra.}
\end{figure}

The twelve irreducible, or {\em pure}, regular apeirohedra fall into a single family, derived from the cubical tessellation in $\E$ and illustrated in the diagram of Figure~\ref{diag} taken from (McMullen \& Schulte, 2002, Ch. 7E). There are a number of relationships between these apeirohedra indicated on the diagram such as {\em duality\/} $\delta$, {\em Petrie duality\/} $\pi$, and {\em facetting\/} $\varphi_2$, as well as the operations $\eta$, $\sigma$ and $\delta\sigma$ which are not further discussed here. The facetting operation $\varphi_2$ applied to a regular polyhedron is reminiscent of the Petrie duality operation, in that the vertices and edges of the polyhedron are retained and the faces are replaced by certain edge-paths, in this case the $2$-holes; here, a {\em $2$-hole}, or simply {\em hole\/}, of a polyhedron is an edge path which leaves a vertex by the second edge from which it entered, always in the same sense (on the left, say, in some local orientation). 

The most prominent apeirohedra of Figure~\ref{diag} are the three Petrie-Coxeter polyhedra $\{4,6\,|\,4\}$, $\{6,4\,|\,4\}$ and $\{6,6\,|\,3\}$, occurring in the top and bottom row; the last entry in the symbols records the length of the holes, $4$ or $3$ in this case, while the first two entries give the standard Schl\"afli symbol. These well-known apeirohedra are the only regular polyhedra in $\E$ with convex faces and skew vertex-figures. 

The Petrie duals of the Petrie-Coxeter polyhedra (related to the Petrie-Coxeter polyhedra under $\pi$) have helical faces given by the Petrie polygons of the original polyhedron. The first subscript in their designation gives the length of their own Petrie polygons, and the second subscript that of their 2-zigzags (edge paths leaving a vertex by the second edge, alternately on the right or left).

Table~\ref{breakdown} is a breakdown of the pure regular apeirohedra by mirror vectors, which also helps understanding why there are exactly 12 examples (McMullen-Schulte, 2002, Ch. 7E). If $\K$ is a regular polyhedron and $R_0,R_1,R_2$ denote the distinguished involutory generators for its symmetry group $G$, the {\em mirror vector\/} $(d_0,d_1,d_2)$ of $\K$ records the dimensions $d_0$, $d_1$ and $d_2$ of the {\em mirrors\/} (fixed point sets) of $R_0$, $R_1$ and $R_2$, respectively. 

It turns out that, mostly due to the irreducibility, only four mirror vectors can occur, namely $(2,1,2)$, $(1,1,2)$, $(1,2,1)$ and $(1,1,1)$. For example, the three apeirohedra with mirror vector $(1,1,1)$ in the last row have a symmetry group generated by three half-turns (reflections in lines) and therefore have only proper isometries as symmetries; these helix-faced apeirohedra occur geometrically in two enantiomorphic forms, yet they are geometrically regular, not chiral!  
 
\begin{table}
\label{breakdown}
\centering
\begin{tabular}{llllll}      
 Mirror    &  $\{3,3\}$ & $\{3,4\}$ & $\{4,3\}$  &faces&vertex-\\
 vector    &                &                &                &         &figures\\
\hline
(2,1,2) & $\{6,6 | 3\}$ & $\{6,4 | 4\}$ & $\{4,6 | 4\}$ &planar & skew  \\
(1,1,2) & $\{\infty,6\}_{4,4}$ & $\{\infty,4\}_{6,4}$ & $\{\infty,6\}_{6,3}$ &helical & skew\\
(1,2,1) & $\{6,6\}_{4}$ & $\{6,4\}_{6}$ & $\{4,6\}_{6}$&  skew &planar \\
(1,1,1) & $\{\infty,3\}^{(a)}$ & $\{\infty,4\}_{\cdot,*3}$ &
$\{\infty,3\}^{(b)}$ &helical &planar  \\
\end{tabular}
\caption{Breakdown of the pure regular apeirohedra by mirror vector.}
\end{table}

While the rows in Table~\ref{breakdown} represent a breakdown of the pure apeirohedra by mirror vector, the first three columns can similarly be seen as grouping the apeirohedra by the crystallographic Platonic solid (the tetrahedron, octahedron or cube) with which each is associated in a manner described below; or, equivalently, as grouping by the corresponding (Platonic) symmetry group. 

To explain this, suppose $\K$ is a pure regular apeirohedron and $G= \langle R_0,R_1,R_2\rangle$, where again $R_0,R_1,R_2$ are the distinguished generators. Then $R_1$ and $R_2$, but not $R_0$, fix the base vertex, the origin $o$ (say), of $\K$. Now consider the translate of the mirror of $R_0$ that passes through $o$, and let $R_0'$ denote the reflection in this translate. If $T$ denotes the translation subgroup of $G$, then $G':=\langle R_0',R_1,R_2\rangle$ is a finite irreducible group of isometries isomorphic to the {\em special group\/} $G\slash T$ of $G$, the quotient of $G$ by its translation subgroup. Now alter (if needed) the generators $R_0',R_1,R_2$ as follows. If a generator is a half-turn (with $1$-dimensional mirror), replace it by the reflection in the plane through $o$ perpendicular to its rotation axis; otherwise leave the generator unchanged. Let $\widehat{R}_{0},\widehat{R}_{1},\widehat{R}_{2}$, respectively, denote the plane reflections derived in this manner from $R_0',R_1,R_2$, and let 
$\widehat{G}:=\langle\widehat{R}_{0},\widehat{R}_{1},\widehat{R}_{2}\rangle$
denote the finite irreducible reflection group in $\E$ generated by them. 

Now since $G$ has to be discrete, $\widehat{G}$ cannot contain $5$-fold rotations. Hence there are only three possibilities for $\widehat{G}$ and its generators, namely $\widehat{G}$ is the symmetry group of the tetrahedron $\{3,3\}$, octahedron $\{3,4\}$ or cube $\{4,3\}$. Bearing in mind that there are just four possible mirror vectors, this then establishes that there are only $12=4\times 3$ pure regular apeirohedra. Note that $G'$ is either again one of these finite reflection groups, or the rotation subgroup of one of these groups (the latter happens only when the mirror vector is $(1,1,1)$).  

Table~\ref{breakdown} also gives details about the geometry of the faces and vertex-figures. It is quite remarkable that in a pure regular apeirohedron with finite faces, the faces and vertex-figures cannot both be planar or both be skew. As we will see, this is very different for chiral polyhedra.

The nets of the pure regular apeirohedra were identified in (O'Keeffe, 2008) and are listed in Table~\ref{purenet}. 

\begin{table}
\label{purenet} 
\centering
\begin{tabular}{lll}
Apeirohedron & &Net \\
\hline
$\{4,6 | 4\}$ & $\{\infty,6\}_{4,4}$ & {\bf pcu} \\
$\{6,4 | 4\}$ & $\{\infty,4\}_{6,4}$ & {\bf sod}\\
$\{6,6 | 3\}$ & $\{\infty,6\}_{6,3}$ & {\bf crs}\\
$\{4,6\}_{6}$ & $\{6,6\}_{4}$ & {\bf hxg} \\
$\{\infty,3\}^{(a)}$ & $\{\infty,3\}^{(b)}$ & {\bf srs} \\
$\{6,4\}_{6}$ & & {\bf nbo}\\
$\{\infty,4\}_{\cdot,*3}$ & & {\bf nbo}\\
\end{tabular}
\caption{The nets of the pure regular apeirohedra. Pairs of Petrie duals share the same net and are listed in the same row. The last two apeirohedra are self-Petrie.}
\end{table}

\section{Chiral Polyhedra}
\label{chipoly}

The classification of chiral polyhedra in ordinary space $\E$ is rather involved. It begins with the observation that chirality, as defined here, does not occur in the classical theory, so in particular there are no convex polyhedra that are chiral.

The chiral polyhedra in $\E$ fall into six infinite families (Schulte, 2004, 2005), each with two or one free parameters  depending on whether the classification is up to congruence or similarity, respectively. Each chiral polyhedron is a non-planar ``irreducible" apeirohedron (with an affinely irrreducible symmetry group), so in particular there are no finite, planar, or  ``reducible" examples. 

The six families comprise three families of apeirohedra with finite skew faces and three families of apeirohedra with infinite helical faces. It is convenient to slightly enlarge each family by allowing the parameters to take certain exceptional values which would make the respective polyhedron regular and, in some cases, finite. The resulting larger families will then contain exactly two regular polyhedra, while all other polyhedra are chiral apeirohedra.

\subsection{Finite-faced polyhedra}

The three families with finite faces only contain apeirohedra and are summarized in Table~\ref{finitefaced}. These apeirohedra are parametrized by two integers which are relatively prime (and not both zero). In fact, the corresponding apeirohedra exist also when the parameters are real, but they are discrete only when the parameters are integers. 

Membership in these families is determined by the basic Schl\"afli symbol, namely $\{6,6\}$, $\{4,6\}$ or $\{6,4\}$. The corresponding apeirohedra are denoted $P(a,b)$, $Q(c,d)$ and $Q(c,d)^*$, respectively, where the star indicates that the apeirohedon $Q(c,d)^*$ in the third family is the dual of the apeirohedron $Q(c,d)$ in the second family. The duality of the  apeirohedra $Q(c,d)$ and $Q(c,d)^*$ is geometric, in that the face centers of one are the vertices of the other. The apeirohedra 
$P(a,b)$ and $P(b,a)$ similarly are geometric duals of each other (again with the roles of vertices and face centers interchanged), and $P(b,a)$ is congruent to $P(a,b)$. Thus $P(a,b)$ is geometrically self-dual, in the sense that its dual $P(a,b)^{*}=P(b,a)$ is congruent to $P(a,b)$. 

\begin{table}
\label{finitefaced}
\centering
\begin{tabular}{llll}      
Schl\"afli &$\{6,6\}$ &$\{4,6\}$  &$\{6,4\}$\\
\hline
Notation & $P(a,b)$& $Q(c,d)$& $Q(c,d)^*$\\
Parameters&$a,b\in\mathbb{Z}$,&$c,d\in\mathbb{Z}$,&$c,d\in\mathbb{Z}$, \\
&$(a,b)=1$ if $b\neq \pm a$&$(c,d)=1$ if $c,d\neq 0$&$(c,d)=1$ if $c,d\neq 0$\\[.04in]
Special goup& $[3,3]^{+}\times \langle -I\rangle$& $[3,4]$&$[3,4]$\\[.04in]
Regular &$P(a,\!-a)\!=\!\{6,\!6\}_{4}$  &$Q(c,\!0)\!=\!\{4,\!6\}_{6}$&$Q(c,\!0)^*\!=\!\{6,\!4\}_{6}$  \\
apeirohedra &$P(a,\!a)\!=\!\{6,\!6|3\}$&$Q(0,\!c)\!=\!\{4,\!6|4\}$&$Q(0,\!c)^*\!=\!\{6,\!4|4\}$  \\[.04in]
&Geom. self-dual&&\\
&$P(a,b)^*\cong P(a,b)$&&\\
\end{tabular}
\caption{The three families of finite-faced chiral apeirohedra and their related regular apeirohedra.}
\end{table}

The chiral apeirohedra in each family have skew faces and skew vertex-figures. However, the two regular apeirohedra in each family have either planar faces or planar vertex-figures. In either case, the distinguished generators $S_1,S_2$ of the symmetry groups $G$ of an apeirohedron in the family are rotatory reflections defined in terms of the parameters $a,b$ or $c,d$, and the symmetry $T:=S_{1}S_{2}$ is a half-turn.

For example, for the apeirohedron $P(a,b)$ of type $\{6,6\}$ the symmetries $S_1$, $S_2$ and $T$ are given by  
\[ \begin{array}{rccl}
S_{1}\colon & (x_{1},x_{2},x_{3}) & \mapsto & (-x_{2},x_{3},x_{1}) + (0,-b,-a),\\ 
S_{2}\colon & (x_{1},x_{2},x_{3}) & \mapsto & -(x_{3},x_{1},x_{2}),  \\
T\colon & (x_{1},x_{2},x_{3}) & \mapsto &  (-x_{1},x_{2},-x_{3}) + (a,0,b).  
\end{array} \]
Then the apeirohedron $P(a,b)$ itself is obtained from the group $G$ generated by $S_1,S_2$ by means of Wythoff's construction as explained above. The base vertex $F_0$ in this case is the origin $o$ (fixed under $S_2$), and the base edge $F_1$ is given by $\{o,u\}$ with 
\[u:= T(o) = (a,0,b).\] 
The base edge $F_1$ lies in the $x_{1}x_{3}$-plane and is perpendicular to the rotation axes of $T$, which in turn is
parallel to the $x_2$-axis and passes through $\frac{1}{2}u$. The base face $F_2$ of $P$ is determined by the orbit of $o$
under $\langle S_1\rangle$ and is given by the generally skew hexagon with vertex-set 
\[ \begin{array}{l}
\{(0,0,0),(0,-b,-a),(b,-a-b,-a),\\
\indent\indent (a+b,-a-b,-a+b),(a+b,-a,b),(a,0,b)\}, 
\end{array}\]
where the vertices are listed in cyclic order. The vertices of $P(a,b)$ adjacent to $o$ are given by the orbit of $u$ under $\langle S_2\rangle$, namely
\[\begin{array}{l}
\{(a,0,b),(-b,-a,0),(0,b,a),\\
\indent\indent (-a,0,-b),(b,a,0),(0,-b,-a)\};
\end{array}\] 
these are the vertices of the generally skew hexagonal vertex-figure of $P(a,b)$ at $o$, listed in the order in which they occur in the vertex-figure. The faces of $P(a,b)$ containing the vertex $o$ are the images of $F_2$ under the elements of $\langle S_2\rangle$. Each face is a generally skew hexagon with vertices given by one half of the vertices of a hexagonal prism. As mentioned earlier, both the faces and vertex-figures are skew if $P(a,b)$ is chiral. 

The apeirohedra $P(a,b)$ are chiral except when $b=\pm a$. If $b=a$ we arrive at the Petrie-Coxeter polyhedron $\{6,6\,|\,3\}$, which has planar (convex) faces but skew vertex-figures. If $b=-a$ we obtain the regular polyhedron $\{6,6\}_4$, which has skew faces but planar (convex) vertex-figures. The vertices of $\{6,6\}_4$ comprise the vertices in one set of alternate vertices of the Petrie-Coxeter polyhedron $\{4,6\,|\,4\}$, while the faces of $\{6,6\}_4$ are the vertex-figures at the vertices in the other set of alternate vertices of $\{4,6\,|\,4\}$.

Table~\ref{finitefaced} also lists the structure of the special groups. Here $[p,q]$ denotes the full symmetry group of a Platonic solid $\{p,q\}$, and $[p,q]^+$ its rotation subgroup; also $-I$ stands for the point reflection in $o$.

\subsection{Helix-faced polyhedra}

The three families of helix-faced chiral apeirohedra and their related regular polyhedra are summarized in Table~\ref{helixfaced}. The corresponding apeirohedra or polyhedra are denoted by $P_1(a,b)$, $P_2(c,d)$ and $P_3(c,d)$. Each family has two {\em real\/}-valued parameters that cannot both be zero. Now the discreteness assumption does not impose any further restrictions on the parameters. 

Membership in these families is determined by the basic Schl\"afli symbol as well as the basic geometry of the helical faces. There are two families of type $\{\infty,3\}$ and one family of type $\{\infty,4\}$. In the first family of type $\{\infty,3\}$ the apeirohedra have helical faces over triangles, and in the second family they have helical faces over squares. Each family contains two regular polyhedra, namely one pure regular apeirohedron as well as one (finite crystallographic) Platonic solid, that is, the tetrahedron, cube or octahedron, respectively,

\begin{table}
\label{helixfaced}
\centering
\begin{tabular}{llll} 
Schl\"afli symbol&$\{\infty,3\}$ &$\{\infty,3\}$&$\{\infty,4\}$\\     
\hline
Notation & $P_1(a,b)$& $P_2(c,d)$& $P_3(c,d)$\\
Parameters&$a,b\in\mathbb{R}, $ &$c,d\in\mathbb{R}, $&$c,d\in\mathbb{R}, $ \\
&$(a,b)\neq (0,0)$ &$(c,d)\neq (0,0)$&$(c,d)\neq (0,0)$\\[.04in]
Helices over& triangles& squares& triangles\\[.04in]
Vertex-figures&triangles&triangles& (planar) squares\\[.04in]
Special group& $[3,3]^+$& $[3,4]^+$&$[3,4]^+$\\[.04in]
Relationships&$P(a,b)^{\varphi_2}$,  & $Q(c,d)^{\varphi_2}$, & $(Q(c,d)^{*})^{\kappa}$ \\
&$a\neq b$&$a\neq 0$ &\\[.02in]
Regular &$P_{1}(a,-a)$&$P_{2}(c,0)$&$P_{3}(0,d)$  \\
polyhedra&$\;\;\;\,=\{\infty,3\}^{(a)}$&$\;\;\;\,=\{\infty,3\}^{(b)}$& $\;\;\;\,=\{\infty,4\}_{\cdot, *3}$\\[.02in]
&$P_{1}(a,a)\!=\!\{3,3\}$&$P_{2}(0,d)\!=\!\{4,3\}$&$P_{3}(c,0)\!=\!\{3,4\}$\\[.05in]
&&&self-Petrie
\end{tabular}
\caption{The three families of helix-faced chiral apeirohedra and their related regular polyhedra.}
\end{table}

The symmetry groups $G$ of the polyhedra $P_1(a,b)$, $P_2(c,d)$ and $P_3(c,d)$ is generated by a screw motion $S_1$ (a rotation followed by a translation along the rotation axis) and an ordinary rotation $S_2$ in an axis through the base vertex $F_{0}:=o$. The screw motion $S_1$ moves the vertices of the helical base face $F_2$ (in an apeirohedron) one step along the face, and $S_2$ applied to the vertex $u$ of the base edge $F_1$ distinct from $o$ produces the planar vertex-figure at $o$. The symmetry $T=S_{1}S_{2}$ is again a half-turn with a rotation axis passing through $\frac{1}{2}u$ and perpendicular to $F_1$.

For example, for the polyhedron $P_2(c,d)$ the symmetries $S_1$, $S_2$ and $T$ are given by 
\[ \begin{array}{rccl}
S_{1}\colon & (x_{1},x_{2},x_{3}) & \mapsto & (-x_{3},x_{2},x_{1}) + (d,c,-c),\\ 
S_{2}\colon & (x_{1},x_{2},x_{3}) & \mapsto & (x_{2},x_{3},x_{1}),  \\
T\colon & (x_{1},x_{2},x_{3}) & \mapsto & (x_{2},x_{1},-x_{3}) + (c,-c,d).  
\end{array} \]
The polyhedron $P_{2}(c,d)$ itself is again obtained from $G=\langle S_1,S_2\rangle$ by Wythoff's construction with base vertex $F_{0}=o$. The base edge $F_1$ is given by $\{0,u\}$ with 
\[ u:=T(o) =  (c,-c,d),\]
and lies in the plane $x_{2}=-x_{1}$. The half-turn $T$ interchanges the vertices $o$ and $u$ of $F_1$, and its rotation axis is parallel (in $\E$) to the line $x_{2}=x_{1}$ in the $x_{1}x_{2}$-plane and perpendicular to the plane $x_{2}=-x_{1}$. The screw motion $S_1$ shifts the base face $F_2$ one step along itself, and since $S_{1}^4$ is the translation along the $x_2$-axis by $(0,4c,0)$, we have helical faces over squares (when $c\neq 0$) spiraling around an axis parallel to the $x_2$-axis. In particular, the vertex-set of $F_2$ is the orbit of $o$ under $\langle S_1\rangle$ and is given by 
\[ \{(c,-c,d),(0,0,0),(d,c,-c),(c+d,2c,-c+d)\} + \mathbb{Z}\!\cdot\! t\]
with $t:=(0,4c,0)$. Here the notation means that the four vectors listed on the left side are successive vertices of $F_2$ whose translates by integral multiples of $t$ comprise all the vertices of~$F_2$.

When $c=0$ we obtain a {\em finite\/} polyhedron $P_{2}(0,d)$, a cube $\{4,3\}$ with a finite group $G$, namely the rotation subgroup of the symmetry group of this cube. In fact, in this case $S_{1}^4$ is the identity mapping and the base face $F_2$ itself is a square (not a helical polygon over a square). 

In the general case, as usual, all other vertices, edges and faces of the apeirohedron (or polyhedron) are the images of $F_{0}$, $F_{1}$ and $F_{2}$ under the group $G$. Moreover, the vertices adjacent to $o$ form the triangular vertex-figure at $o$ and are given by the three points 
\[ u = (c,-c,d),\, S_{2}(u)=(-c,d,c),\, S_{2}^2(u)=(d,c,-c)\] 
in the plane $x_{1}+x_{2}+x_{3}=d$. 

\subsection{General properties}

The six families of chiral apeirohedra and related regular polyhedra have stunning geometric and combinatorial properties and exhibit some rather unexpected phenomena. 

For all three families of apeirohedra with finite faces it is almost true that different parameter values give combinatorially non-isomorphic apeirohedra (that is, the underlying abstract apeirohedra are non-isomorphic). More explicitly, $P(a,b)$ and $P(a',b')$ are combinatorially isomorphic if and only if 
\[ (a',b')=\pm (a,b),\pm (b,a); \] 
and similarly, $Q(c,d)$ and $Q(c',d')$, and  hence $Q(c,d)^*$ and $Q(c',d')^*$, are combinatorially isomorphic if and only if 
\[ (c',d')=\pm (c,d),\pm (-c,d).\] 

This phenomenon is perhaps even more surprising when expressed in terms of the similarity classes of the apeirohedra within each family, which are parametrized by a single rational parameter, namely $a/b$ or $c/d$ (when $b,d\neq 0$). These similarity classes exhibit a very strong discontinuity, in that any small positive change in the rational parameter produces a new similarity class in which the apeirohedra are not combinatorially isomorphic to those in the original similarity class.

By contrast, in the three families of helix-faced apeirohedra and related polyhedra, each chiral apeirohedron is combinatorially isomorphic to the regular apeirohedron in its family, so in particular it is combinatorially regular, but not geometrically regular (Pellicer \& Weiss, 2010). In fact, the chiral apeirohedra in each family can be viewed as ``chiral deformations" of the regular apeirohedron in this family. At the other extreme they also allow a ``deformation" to the finite regular polyhedron (Platonic solid) in the family. 

On the other hand, the finite-faced chiral apeirohedra are also combinatorially chiral (Pellicer \& Weiss, 2010). In other words, these apeirohedra are intrinsically chiral and thus not combinatorially isomorphic to a regular apeirohedron in their family.

In some sense, each chiral helix-faced apeirohedron can be thought of as unraveling the Platonic solid in its family. In fact, for all parameter values $a,b$ and $c,d$ we have the following coverings of polyhedra:
\[ P_{1}(a,b)\! \mapsto\! \{3,3\},\; P_{2}(c,d)\! \mapsto\! \{4,3\},\; P_{3}(c,d)\! \mapsto\! \{3,4\}. \]
Informally speaking, under these coverings each helical face is ``compressed" (like a spring) to become a polygon (triangle or square) over which it has been rising.

The nets arising as edge graphs of chiral polyhedra in $\E$ will be analyzed in a forthcoming paper.

\section{Regular Polygonal Complexes}
\label{regpolycom}

We now turn to polygonal complexes with possibly more than two faces meeting at an edge. All regular polygonal complexes that are not polyhedra turn out to be infinite and have an affinely irreducible symmetry group. 

Unlike a regular polyhedron, a regular polygonal complex $\K$ can have a symmetry group that is transitive, but not simply transitive, on the flags. As we mentioned earlier, $\K$ then has face mirrors, meaning that $\K$ has planar faces and that each flag stabilizer is generated by the reflection in the plane through the (planar) face in the flag. 

The classification of regular polygonal complexes naturally breaks down into two cases, namely the enumeration of the simply flag-transitive complexes and that of the non-simply flag-transitive complexes (Pellicer \& Schulte, 2010, 2013). 

All regular polyhedra, finite or infinite, are simply flag-transitive polygonal complexes.

\subsection{Non-simply flag-transitive complexes as 2-skeletons}
\label{nonsim}

For a regular polygonal complex $\K$ with a non-simply flag transitive symmetry group $G$, the existence of face mirrors allows us to recognize $\K$ as the $2$-skeleton of a certain type of incidence structure of rank $4$ in $\E$, called a {\em $4$-apeirotope} (McMullen \& Schulte, 2002, Ch. 7F). Thus $\K$ consists of the vertices, edges, and faces (of rank 2) of this $4$-apeirotope. The $2$-skeleton of the cubical tessellation shown in Figure~\ref{skel2} is an example of a regular polygonal complex of this kind, and the underlying cubical tessellation is the corresponding 4-apeirotope. 

The $4$-apeirotopes involved are themselves {\em regular\/}, in the sense that they have a flag-transitive symmetry group on their own (coinciding with $G$); in fact, the generating reflection of the base flag stabilizer for $\K$ is the fourth involutory generator needed to suitably generate the symmetry group of this rank 4 structure. 

There are eight regular $4$-apeirotopes in $\E$, occurring in four pairs of ``Petrie-duals" (McMullen \& Schulte, 2002, Ch. 7F). The apeirotopes in each pair have the same 2-skeleton and produce the same polygonal complex. Thus up to similarity there are four non-simply flag-transitive complexes in $\E$, each with the same symmetry group as its two respective apeirotopes. The $2$-skeleton of the cubical tessellation has square faces and is the only example of a non-simply flag-transitive regular polygonal complex with finite faces. The three other complexes all have (planar) zigzag faces, with either $3$ or $4$ faces meeting at each edge.

\begin{table}
\label{skelcom}
\centering
\begin{tabular}{llllllll}
Apeirotope&&$\!\!r$ &$\!\!\rm{Vertex}$&Net \\
&&&$\!\! \rm{-Figure}$&\\
\hline
$\{4, 3, 4\}$ &$\!\!\! \{\{4, 6 \mid 4\}, \{6, 4\}_3\}$ &$\!\! 4$&$\!\! \rm{octahedron}$& {\bf pcu}\\
$\{\{\infty, 3\}_6 \# \{ \, \}, \{3, 3\}\}$
&$\!\!\! \{\{\infty, 4\}_4 \# \{\infty\}, \{4, 3\}_3\}$  &$\!\! 3$&$\!\! \rm{tetrahedron}$& {\bf dia} \\
$\{\{\infty, 3\}_6 \# \{ \, \}, \{3, 4\}\}$
&$\!\!\! \{\{\infty, 6\}_3 \# \{\infty\}, \{6, 4\}_3\}$ &$\!\! 4$&$\!\! \rm{octahedron}$& {\bf pcu}\\
$\{\{\infty, 4\}_4 \# \{\, \}, \{4, 3\}\}$ 
&$\!\!\! \{\{\infty, 6\}_3 \# \{\infty\}, \{6, 3\}_4\}$ &$\!\! 3$&$\!\! \rm{cube}$& {\bf bcu}
\end{tabular}
\caption{The four nets arising as edge graphs of regular $4$-apeirotopes in $\mathbb{E}^3$. Pairs of Petrie-dual apeirotopes are listed in the same row. }
\end{table}

In Table~\ref{skelcom} we list the four pairs of Petrie dual regular $4$-apeirotopes in $\mathbb{E}^3$, along with  information about the regular polygonal complexes $\K$ arising as their 2-skeletons. In particular we give the number $r$ of faces at an edge, as well as the structure of the vertex-figure of $\K$; here an entry in the vertex-figure column of Table~\ref{skelcom} listing a Platonic solid is meant to represent the (geometric) edge-graph of this solid. Note that in each case the vertex-figure of the polygonal complex $\K$ is simply the edge graph of the Platonic solid which occurs as the vertex-figure of the regular $4$-apeirotope listed in the first column. Thus the vertex-figure column gives the structure of the vertex-figures of both the regular polygonal complex and the corresponding regular $4$-apeirotope in the first column. Similarly, the number $r$ for $\K$ also coincides with the number of faces at an edge of the corresponding $4$-apeirotope. 

Table~\ref{skelcom} also identifies the underlying net of each complex, which here coincides with the edge graphs of the two related regular $4$-apeirotopes.

\subsection{Simply flag-transitive complexes}
\label{simply}

In addition to the 48 regular polyhedra described in Section~\ref{regpoly} there are $21$ simply flag-transitive regular polygonal complexes in $\mathbb{E}^3$ which are not polyhedra, up to similarity. This gives a total of $69$ simply flag-transitive regular complexes, up to similarity and relative scaling of components for blended polyhedra.

Let $\K$ be a simply flag-transitive complex which is not a polyhedron, that is, $\K$ has $r\geq 3$ faces at each edge. Then the generating subgroups $G_0,G_1,G_2$ of its symmetry group $G$ described earlier take a very specific form:\ $G_0$ is generated by a point, line, or plane reflection $R_0$; the subgroup $G_1$ is generated by a line or plane reflection $R_1$; and $G_2$ is a cyclic or dihedral group of order $r$. The {\em mirror vector\/} $(d_0,d_1)$ of $\K$ records the dimensions $d_0$ and $d_1$ of the mirrors of $R_0$ and $R_1$, respectively. (Note that, for a polyhedron, $G_2$ would also be generated by a reflection and the complete mirror vector introduced earlier would list the dimensions of all three mirrors.)

Table~\ref{tabsimply} summarizes the $21$ simply flag-transitive complexes and some of their properties. In writing $\K_i(j,k)$ for a complex, $(j,k)$ indicates its mirror vector and $i$ is its label (serial number) in the list of regular complexes with the same mirror vector $(j,k)$. There are columns for the pointwise edge stabilizer $G_2$, the number $r$ of faces at each edge, the types of faces and vertex-figures, the vertex-set, the special group $G^*$, and the corresponding net. In the face column we use symbols like $p_c$, $p_s$, $\infty_2$, or $\infty_k$ with $k=3$ or $4$, respectively, to indicate that the faces are \underbar{c}onvex $p$-gons, \underbar{s}kew $p$-gons, planar zigzags, or helical polygons over $k$-gons. (A planar zigzag is viewed as a helix over a ``$2$-gon", where here a $2$-gon is a line segment traversed in both directions. Hence the use of $\infty_2$.) An entry in the vertex-figure column describing a solid figure in $\E$ is meant to represent the geometric edge-graph of this figure, with ``double" indicating the double edge-graph (the edges have multiplicity 2). The entry ``ns-cuboctahedron" in the vertex-figure column stands for the edge graph of a ``non-standard cuboctahedron", a certain realization with skew square faces of the ordinary cuboctahedron.

\begin{table}
\label{tabsimply}
\centering
\begin{tabular}{llllllll}
Complex&$G_2$ & $r$ &Face&Vertex&$\!\rm{Vertex}$& $\!G^*$&$\! \rm{Net}$\\
              &            &       &        &-Figure&-Set&&\\
\hline
$\K_1(1,2)$& $D_2$  & $4$ &$4_s$ & cuboctahedron&$\Lambda_{2}$&$\![3,4]$&$\!{\bf fcu}$\\
$\K_2(1,2)$& $C_3$& $3$ &$4_s$ & cube&$\Lambda_{3}$&$\![3,4]$&$\!{\bf bcu}$\\
$\K_3(1,2)$& $D_3$& $6$   &$4_s$ &double cube&$\Lambda_{3}$&$\![3,4]$&$\!{\bf bcu}$\\
$\K_4(1,2)$& $D_2$& $4$ & $6_s$&octahedron&$\Lambda_{1}$&$\![3,4]$&$\!{\bf pcu}$\\
$\K_5(1,2)$& $D_2$& $4$ &$6_s$&double square&$V$&$\![3,4]$&$\!{\bf nbo}$\\
$\K_6(1,2)$& $D_4$& $8$ &$6_s$&double octahedron&$\Lambda_{1}$&$\![3,4]$&$\!{\bf pcu}$\\
$\K_7(1,2)$& $D_3$& $6$ &$6_s$&double tetrahedron&$W$&$\![3,4]$&$\!{\bf dia}$\\
$\K_8(1,2)$& $D_2$& $4$ &$6_s$&cuboctahedron&$\Lambda_{2}$&$\![3,4]$&$\!{\bf fcu}$\\[.07in]
%\hline
$\K_1(1,1)$& $D_3$& $6$ &$\infty_3$&double cube&$\Lambda_{3}$&$\![3,4]$&$\!{\bf bcu}$\\
$\K_2(1,1)$& $D_2$& $4$ &$\infty_3$&double square&$V$&$\![3,4]$&$\!{\bf nbo}$\\
$\K_3(1,1)$& $D_4$& $8$ &$\infty_3$&double octahedron&$\Lambda_{1}$ & $\![3,4]$&$\!{\bf pcu}$ \\
$\K_4(1,1)$& $D_3$ & $6$ &$\infty_4$& double tetrahedron&$W$ &$\![3,4]$&$\!{\bf dia}$\\
$\K_5(1,1)$& $D_2$& $4$ &$\infty_4$&ns-cuboctahedron&$\Lambda_{2}$&$\![3,4]$&$\!{\bf fcu}$ \\
$\K_6(1,1)$& $C_3$& $3$ &$\infty_4$&tetrahedron& $W$&$\![3,4]^+$&$\!{\bf dia}$\\
$\K_7(1,1)$& $C_4$& $4$ &$\infty_3$&octahedron&$\Lambda_{1}$ & $\![3,4]^+$&$\!{\bf pcu}$\\
$\K_8(1,1)$& $D_2$& $4$ &$\infty_3$&ns-cuboctahedron&$\Lambda_{2}$&$\![3,4]$&$\!{\bf fcu}$ \\
$\K_9(1,1)$& $C_3$& $3$ &$\infty_3$&cube&$\Lambda_{3}$&$\![3,4]^+$&$\!{\bf bcu}$\\[.07in]
%\hline
$\K(0,1)$& $D_2$& $4$ &$\infty_2$&ns-cuboctahedron&$\Lambda_{2}$&$\![3,4]$&$\!{\bf fcu}$ \\
%\hline
$\K(0,2)$& $D_2$& $4$ &$\infty_2$&cuboctahedron& $\Lambda_{2}$&$\![3,4]$&$\!{\bf fcu}$ \\
%\hline
$\K(2,1)$& $D_2$& $4$ & $6_c$ &ns-cuboctahedron&$\Lambda_{2}$&$\![3,4]$ &$\!{\bf fcu}$\\
%\hline
$\K(2,2)$& $D_2$& $4$ &$3_c$&cuboctahedron&$\Lambda_{2}$&$\![3,4]$&$\!{\bf fcu}$
\end{tabular}
\caption{The 21 simply flag-transitive regular complexes in $\mathbb{E}^3$ which are not polyhedra, and their nets.}
\end{table}

For all but three complexes, the special group $G^*$ is the full octahedral group $[3,4]$. In the three exceptional cases $G^*$ is the octahedral rotation group $[3,4]^+$; note that in these cases we must have a mirror vector $(1,1)$ and a cyclic group $G_{2}$.

For all but five complexes of Table~\ref{tabsimply} the vertex-set is a lattice, namely one of the following, up to scaling:\ the standard (``primitive") {\em cubic\/} lattice $\Lambda_{1}:=\mathbb{Z}^{3}$; the {\em face-centered cubic\/} lattice $\Lambda_{2}$, with basis $(1,1,0)$, $(-1,1,0)$, $(0,-1,1)$, consisting of all integral vectors with even coordinate sum; or the {\em body-centered cubic\/} lattice $\Lambda_3$, with basis $(2,0,0)$, $(0,2,0)$, $(1,1,1)$. In the five exceptional cases the vertex-set is either  
\[ V:=\Lambda_{1}\!\setminus\! ((0,0,1)\!+\!\Lambda_{3})\]
or
\[ W:= 2\Lambda_{2} \cup ((1,-1,1)\!+\!2\Lambda_{2}), \]
again up to scaling.

The last column of Table~\ref{tabsimply} lists the nets of the simply flag-transitive regular complexes which are not polyhedra. The coordination figures of the nets are just the convex hulls of the vertex-figures of the corresponding regular complex, ignoring multiplicity of the edges in the vertex-figure if it occurs. For example, the vertex-figure of the complex $\mathcal{K}_{1}(1,2)$ is the graph of an ordinary cuboctahedron and so the coordination figure of the net for $\mathcal{K}_{1}(1,2)$ is a cuboctahedron; hence the net is quasiregular and must coincide with the face-centered cubic lattice $\bf fcu$, which we denoted by $\Lambda_2$. The nets ${\bf pcu}$ and ${\bf bcu}$ similarly correspond to the lattices $\Lambda_1$ and $\Lambda_3$, respectively. 

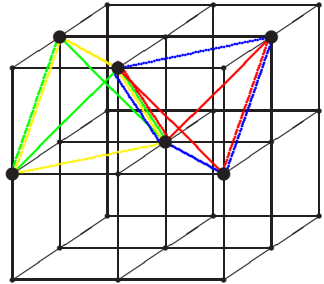
\begin{figure}
\centering
\begin{picture}(130,110)
\multiput(10,0)(0,40){3}{
\begin{picture}(100,25)
\thinlines
\multiput(0,0)(40,0){3}{\circle*{2}}
\multiput(18,12)(40,0){3}{\circle*{2}}
\multiput(36,24)(40,0){3}{\circle*{2}}
\multiput(0,0)(18,12){3}{\line(1,0){80}}
\multiput(0,0)(40,0){3}{\line(3,2){36}}
\end{picture}}
\put(10,0){
\begin{picture}(180,140)
\thinlines
\multiput(0,0)(40,0){3}{\line(0,1){80}}
\multiput(18,12)(40,0){3}{\line(0,1){80}}
\multiput(36,24)(40,0){3}{\line(0,1){80}}
\end{picture}}
\put(10,0){
\begin{picture}(100,80)
\thicklines
\put(0,40){\color{green}\line(1,1){40}}
\multiput(40,81)(.45,-.7){40}{\color{green}\circle*{1}}          
\put(58,52){\color{green}\line(-1,1){40}}
\multiput(18,93.5)(-.3,-.866){60}{\color{green}\circle*{1}}
\end{picture}}
\put(10,0){
\begin{picture}(100,80)
\thicklines
\multiput(0,40)(0.98,.2){60}{\color{yellow}\circle*{1}} 
\multiput(40,79)(.45,-.7){40}{\color{yellow}\circle*{1}} 
\multiput(40,80)(-.55,.3){40}{\color{yellow}\circle*{1}}         
\multiput(18,90)(-.3,-.866){60}{\color{yellow}\circle*{1}}
\end{picture}}
\put(10,0){
\begin{picture}(100,80)
\thicklines
\put(40,80){\color{red}\line(1,-1){40}}
\multiput(80,42)(.3,.866){60}{\color{red}\circle*{1}}
\put(58,52){\color{red}\line(1,1){40}}
\multiput(40,82.5)(.45,-.68){42}{\color{red}\circle*{1}}  
\end{picture}}
\put(10,0){
\begin{picture}(100,80)
\thicklines
\multiput(80,40)(-.55,.3){40}{\color{blue}\circle*{1}}      
\multiput(80,38)(.3,.866){60}{\color{blue}\circle*{1}}
\multiput(40,80)(0.98,.2){60}{\color{blue}\circle*{1}} 
\multiput(39,77.5)(.45,-.68){40}{\color{blue}\circle*{1}}  
\end{picture}}
\put(10,0){
\begin{picture}(80,80)
\put(0,40){\circle*{5}}
\put(18,92){\circle*{5}}
\put(80,40){\circle*{5}}
\put(98,92){\circle*{5}}
\put(40,80){\circle*{5}}
\put(58,52){\circle*{5}}
\end{picture}}
\end{picture}
\caption{The four skew square faces of $\mathcal{K}_{1}(1,2)$ sharing an edge. Each face is a Petrie polygon of a regular tetrahedron inscribed in a cube of the cubical tessellation. The tetrahedra in adjacent cubes have different orientations. The net is {\bf fcu}.}
\label{figk112}
\end{figure}

Figure~\ref{figk112} illustrates a local picture of $\mathcal{K}_{1}(1,2)$ around an edge, showing the four skew square faces meeting at the edge. The entire complex can be thought of as being built from an infinite family of Petrie duals of regular tetrahedra inscribed in cubes, one per cube of the cubical tessellation, such that the Petrie duals in adjacent cubes share a common edge and have opposite orientation. Thus the complex has infinitely many ``small" finite subcomplexes each a regular polyhedron in itself.  

\begin{figure}
\centering
\begin{picture}(140,110)
\multiput(10,0)(0,40){3}{
\begin{picture}(100,25)
\thinlines
\multiput(0,0)(40,0){3}{\circle*{3}}
\multiput(18,12)(40,0){3}{\circle*{3}}
\multiput(36,24)(40,0){3}{\circle*{3}}
\multiput(0,0)(18,12){3}{\line(1,0){80}}
\multiput(0,0)(40,0){3}{\line(3,2){36}}
\end{picture}}
\put(10,0){
\begin{picture}(180,140)
\thinlines
\multiput(0,0)(40,0){3}{\line(0,1){80}}
\multiput(18,12)(40,0){3}{\line(0,1){80}}
\multiput(36,24)(40,0){3}{\line(0,1){80}}
\end{picture}}
\put(10,0){
\begin{picture}(100,80)
\thicklines
\put(41,40){\color{green}\line(0,1){40}}
\put(40,81){\color{green}\line(1,0){40}}
\multiput(80,81)(.45,.3){40}{\color{green}\circle*{1}}          
\put(97,52){\color{green}\line(0,1){40}}
\put(98,51){\color{green}\line(-1,0){40}}
\multiput(58,51)(-.45,-.3){40}{\color{green}\circle*{1}}
\end{picture}}
\put(10,0){
\begin{picture}(100,80)
\thicklines
\put(58,54){\color{red}\line(1,0){40}}
\put(99,52){\color{red}\line(0,1){40}}
\multiput(98,93)(.45,.3){40}{\color{red}\circle*{1}}    
\put(116,105){\color{red}\line(-1,0){40}}
\put(77,104){\color{red}\line(0,-1){40}}
\multiput(58,54.2)(.45,.3){40}{\color{red}\circle*{1}}
\put(58,50){\color{blue}\line(-1,0){40}} 
\put(19,52){\color{blue}\line(0,1){40}}
\multiput(18,93)(.45,.3){40}{\color{blue}\circle*{1}}
\put(36,105){\color{blue}\line(1,0){40}}
\put(75,104){\color{blue}\line(0,-1){40}}
\multiput(76,61.6)(-.45,-.3){40}{\color{blue}\circle*{1}}
\multiput(58,51)(.45,.3){40}{\color{green}\circle*{1}}
\put(75,64){\color{green}\line(0,-1){40}}
\put(76,23){\color{green}\line(-1,0){40}}
\multiput(36,23)(-.45,-.3){40}{\color{green}\circle*{1}}
\put(17,51){\color{green}\line(1,0){40}}
\put(19,52){\color{green}\line(0,-1){40}}
\multiput(58,53)(-.45,-.3){40}{\color{yellow}\circle*{1}}
\put(58,53){\color{yellow}\line(-1,0){40}}
\put(39,40){\color{yellow}\line(0,1){40}}
\put(40,81){\color{yellow}\line(-1,0){40}}
\multiput(0,81)(.45,.3){40}{\color{yellow}\circle*{1}}
\put(17,52){\color{yellow}\line(0,1){40}}
\put(58,54){\color{red}\line(-1,0){40}}
\put(17,52){\color{red}\line(0,-1){40}}
\multiput(18,11)(-.45,-.3){40}{\color{red}\circle*{1}}    
\put(0,-1){\color{red}\line(1,0){40}}
\put(39,0){\color{red}\line(0,1){40}}
\multiput(40,42)(.45,.3){40}{\color{red}\circle*{1}}
\put(58,50){\color{blue}\line(1,0){40}} 
\put(97,52){\color{blue}\line(0,-1){40}}
\multiput(98,11)(-.45,-.3){40}{\color{blue}\circle*{1}}
\put(80,-1){\color{blue}\line(-1,0){40}}
\put(41,0){\color{blue}\line(0,1){40}}
\multiput(40,38)(.45,.3){40}{\color{blue}\circle*{1}}
\multiput(58,53)(.45,.3){40}{\color{yellow}\circle*{1}}
\put(77,64){\color{yellow}\line(0,-1){40}}
\put(76,23){\color{yellow}\line(1,0){40}}
\put(99,12){\color{yellow}\line(0,1){40}}
\multiput(116,24)(-.45,-.3){40}{\color{yellow}\circle*{1}}
\put(58,53){\color{yellow}\line(1,0){40}}
\put(58,52){\circle*{6}}
\put(18,52){\circle*{6}}
\put(98,52){\circle*{6}}
\put(40,40){\circle*{6}}
\put(76,64){\circle*{6}}
\end{picture}}
\multiput(10,0)(0,40){3}{
\begin{picture}(100,25)
\thinlines
\multiput(0,0)(40,0){3}{\circle*{3}}
\multiput(18,12)(40,0){3}{\circle*{3}}
\multiput(36,24)(40,0){3}{\circle*{3}}
\multiput(0,0)(18,12){3}{\line(1,0){80}}
\multiput(0,0)(40,0){3}{\line(3,2){36}}
\end{picture}}

\end{picture}
\caption{The faces of $\mathcal{K}_{5}(1,2)$ are Petrie polygons of cubes in the cubical tessellation of $\E$. For every cube of the tessellation, only one of its Petrie polygons is a face of $\mathcal{K}_{5}(1,2)$. Shown are the eight faces of $\mathcal{K}_{5}(1,2)$ that have a vertex in common, here located at the center. Each edge containing this vertex lies in four faces. The vertex-figure of $\K_{5}(1,2)$ at the central vertex is the double square spanned by the four outer black nodes in the horizontal plane through the center. The net is {\bf nbo}.}
\label{figk512}
\end{figure}
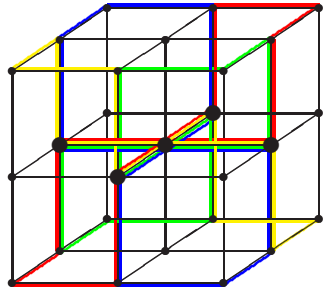

Figures~\ref{figk412} and \ref{figk512} depict the local structure around a vertex for two further simply flag-transitive regular complexes with mirror vector $(1,2)$, namely $\K_{4}(1,2)$ and $\K_{5}(1,2)$. They also are related to the cubical tessellation and their nets are {\bf pcu} and {\bf nbo}, respectively.

It turns out that the edge graph of each simply-flag-transitive regular complex $\K$ is a regular or quasiregular net in the sense described earlier. This can be seen as follows. It is clear that the net is uninodal, since the symmetry group of $\K$ acts transitively on the vertices of $\K$ and is a subgroup of the symmetry group of the net.

If the vertex-figure of $\K$ is a cuboctahedron or non-standard (ns) cuboctahedron, then the coordination figures of the net are convex cuboctahedra. Hence the net is quasiregular and must coincide with {\bf fcu}. Note here that the argument also applies if the vertex-figures of $\K$ are ns-cuboctahedra, that is, when $\K$ is one of the four complexes $\mathcal{K}_5(1,1)$, $\mathcal{K}_8(1,1)$, $\mathcal{K}(0,1)$ and $\mathcal{K}(2,1)$. In these four cases, the edges in the vertex-figure connect pairs of vertices of the standard cuboctahedron which are midpoints of edges two steps apart on a Petrie polygon of the cube used to construct the cuboctahedron as the convex hull of the midpoints of its edges. However, an edge in the vertex-figure of a polygonal complex at a given vertex represents a face of the complex that contains this vertex, in that it joins the two vertices of the face that are adjacent to the given vertex in the complex. Thus the edges of the vertex-figure capture the faces of the complex, not its underlying net. The local structure of the net is only determined by the vertices of the vertex-figure, not its edges.

For all other regular complexes $\K$, except $\K_2(1,2)$ and $\K_{4}(1,2)$, we can appeal to the classification of regular nets. In fact, the coordination figures of the net of $\K$ are easily seen to be squares, tetrahedra, octahedra or cubes, and the vertex-figure subgroup $G_{F_0}=\langle R_{1},G_2\rangle$ of the symmetry group $G$ at the base vertex of $\K$ contains (at least) the rotation symmetry group of the coordination figure at this vertex of the net. In the case of the square, the rotational symmetries are taken relative to $\E$. Thus the complex $\K$ has enough symmetries to guarantee that its edge graph is a regular net.

The situation changes for the complexes $\K_2(1,2)$ and $\K_{4}(1,2)$ with cubes and octahedra as coordination figures, respectively. Now the vertex-figure subgroups are too small to imply regularity of the net based on symmetries of $\K$; these subgroups are given by $[3,3]^{+}\times \langle -I\rangle$ and $[3,3]$, respectively. However, $\K_{2}(1,2)$ is a subcomplex of $\K_{3}(1,2)$ with the same edge graph, and the latter is the regular net {\bf bcu} by our previous analysis. Thus the net of $\K_2(1,2)$ also is {\bf bcu}. Similarly, $\K_{4}(1,2)$ is a subcomplex of $\K_{6}(1,2)$ with the same edge graph, namely the edge graph of the cubical tessellation of $\E$ (see Figure~\ref{figk412}), which is the regular net {\bf pcu}.

Thus the edge graphs of all regular polygonal complexes which are not polyhedra, simply flag-transitive or not, are regular or quasiregular nets. The regular net {\bf srs} does not occur in this but all others do; however, {\bf srs} is the net of two pure regular polyhedra with helical faces (see Table~\ref{purenet}).  

\begin{figure}
\centering
\begin{picture}(140,110)
\multiput(10,0)(0,40){3}{
\begin{picture}(100,25)
\thinlines
\multiput(0,0)(40,0){3}{\circle*{3}}
\multiput(18,12)(40,0){3}{\circle*{3}}
\multiput(36,24)(40,0){3}{\circle*{3}}
\multiput(0,0)(18,12){3}{\line(1,0){80}}
\multiput(0,0)(40,0){3}{\line(3,2){36}}
\end{picture}}
\put(10,0){
\begin{picture}(180,140)
\thinlines
\multiput(0,0)(40,0){3}{\line(0,1){80}}
\multiput(18,12)(40,0){3}{\line(0,1){80}}
\multiput(36,24)(40,0){3}{\line(0,1){80}}
\end{picture}}
\put(10,0){
\begin{picture}(100,80)
\thicklines
\put(40,-1){\color{green}\line(1,0){40}}
\multiput(80,0)(.45,.3){40}{\color{green}\circle*{1}}          
\put(98,12){\color{green}\line(0,1){40}}
\put(98,51){\color{green}\line(-1,0){40}}
\multiput(58,50.5)(-.45,-.3){40}{\color{green}\circle*{1}}
\put(40,40){\color{green}\line(0,-1){40}}
\end{picture}}
\put(10,0){
\begin{picture}(100,80)
\thicklines
\put(58,49){\color{blue}\line(1,0){40}} 
\multiput(98,52)(.45,.3){40}{\color{blue}\circle*{1}}          
\put(116,64){\color{blue}\line(0,1){40}}
\put(116,104){\color{blue}\line(-1,0){40}}
\multiput(76,104)(-.45,-.3){40}{\color{blue}\circle*{1}}
\put(61,92){\color{blue}\line(0,-1){40}}
\end{picture}}
\put(10,0){
\begin{picture}(100,80)
\thicklines
\put(40,1){\color{yellow}\line(1,0){40}}
\put(80,0){\color{yellow}\line(0,1){40}}
\multiput(80,40)(.45,.3){40}{\color{yellow}\circle*{1}}    
\put(98,53.5){\color{yellow}\line(-1,0){40}}
\put(56.5,52){\color{yellow}\line(0,-1){40}}
\multiput(40,0)(.45,.3){40}{\color{yellow}\circle*{1}}
\end{picture}}
\put(10,0){
\begin{picture}(100,80)
\thicklines
\put(58,55){\color{red}\line(1,0){40}}
\put(98,52){\color{red}\line(0,1){40}}
\multiput(98,92)(.45,.3){40}{\color{red}\circle*{1}}    
\put(116,105){\color{red}\line(-1,0){40}}
\put(76,104){\color{red}\line(0,-1){40}}
\multiput(58,55)(.45,.3){40}{\color{red}\circle*{1}}
\end{picture}}
\put(10,0){
\begin{picture}(100,80)
\thicklines
\put(56.5,52){\color{yellow}\line(0,1){40}}
\multiput(58,53.5)(.45,.3){40}{\color{yellow}\circle*{1}}  
\put(58,92){\color{yellow}\line(1,0){40}}
\multiput(98,94)(.45,.3){40}{\color{yellow}\circle*{1}}
\put(76,64){\color{yellow}\line(1,0){40}}
\put(117,64){\color{yellow}\line(0,1){40}}
\put(61,52){\color{blue}\line(0,-1){40}}
\multiput(58,49)(-.45,-.3){40}{\color{blue}\circle*{1}}
\put(40,40){\color{blue}\line(1,0){40}}
\put(79,40){\color{blue}\line(0,-1){40}}
\multiput(80,-1)(.45,.3){40}{\color{blue}\circle*{1}}  
\put(98,12){\color{blue}\line(-1,0){40}}
\put(58,55){\color{red}\line(-1,0){40}} 
\multiput(18,51)(.45,.3){40}{\color{red}\circle*{1}}
\multiput(58,13.5)(.45,.3){40}{\color{red}\circle*{1}}
\put(36,23){\color{red}\line(1,0){40}}
\put(37,24){\color{red}\line(0,1){40}}
\put(55,52){\color{red}\line(0,1){40}}
\put(58,49){\color{blue}\line(-1,0){40}} 
\put(18,12){\color{blue}\line(0,1){40}}
\multiput(58,49)(.45,.3){40}{\color{blue}\circle*{1}}
\put(76,64){\color{blue}\line(0,-1){40}}
\put(76,25){\color{blue}\line(-1,0){40}}
\multiput(36,24)(-.45,-.3){40}{\color{blue}\circle*{1}}
\multiput(58,50.5)(.45,.3){40}{\color{green}\circle*{1}}
\multiput(18,13)(.45,.3){40}{\color{green}\circle*{1}}
\put(18,12){\color{green}\line(1,0){40}}
\put(59.5,12){\color{green}\line(0,1){40}}
\put(76,64){\color{green}\line(-1,0){40}}
\put(36,64){\color{green}\line(0,-1){40}}
\multiput(58,53.5)(-.45,-.3){40}{\color{yellow}\circle*{1}}
\put(58,53.5){\color{yellow}\line(-1,0){40}}
\put(40,40){\color{yellow}\line(0,1){40}}
\put(40,80){\color{yellow}\line(-1,0){40}}
\multiput(0,80)(.45,.3){40}{\color{yellow}\circle*{1}}
\put(18,92){\color{yellow}\line(0,-1){40}}
\put(40,40){\color{red}\line(-1,0){40}}
\put(1,40){\color{red}\line(0,1){40}}
\multiput(40,42)(.45,.3){40}{\color{red}\circle*{1}}
\put(58,92){\color{red}\line(-1,0){40}}
\multiput(18,93)(-.45,-.3){40}{\color{red}\circle*{1}}
\put(55,53){\color{red}\line(0,-1){40}}
\put(0,39){\color{green}\line(0,1){40}}
\put(0,81){\color{green}\line(1,0){40}}
\multiput(40,80)(.45,.3){40}{\color{green}\circle*{1}}
\multiput(0,40)(.45,.3){40}{\color{green}\circle*{1}}
\put(59.5,92){\color{green}\line(0,-1){40}}
\put(58,50.5){\color{green}\line(-1,0){40}}
\put(58,52){\circle*{6}}
\put(18,52){\circle*{6}}
\put(98,52){\circle*{6}}
\put(58,12){\circle*{6}}
\put(58,92){\circle*{6}}
\put(40,40){\circle*{6}}
\put(76,64){\circle*{6}}
\end{picture}}
\multiput(10,0)(0,40){3}{
\begin{picture}(100,25)
\thinlines
\multiput(0,0)(40,0){3}{\circle*{3}}
\multiput(18,12)(40,0){3}{\circle*{3}}
\multiput(36,24)(40,0){3}{\circle*{3}}
\multiput(0,0)(18,12){3}{\line(1,0){80}}
\multiput(0,0)(40,0){3}{\line(3,2){36}}
\end{picture}}
\end{picture}
\caption{The faces of $\mathcal{K}_{4}(1,2)$ are the Petrie polygons of alternate cubes in the cubical tessellation of $\E$. For every cube occupied, all its Petrie polygons occur as faces of $\mathcal{K}_{4}(1,2)$. Shown are the twelve faces of $\mathcal{K}_{4}(1,2)$ that have a vertex in common, here located at the center. Each edge containing this vertex lies in four faces. The vertex-figure of $\K_{4}(1,2)$ at the central vertex is the octahedron spanned by the six outer black nodes. The net is {\bf pcu}.}
\label{figk412}
\end{figure}
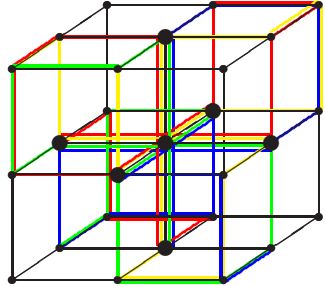

\noindent
{\sl Acknowledgements.}\
I am grateful to Davide Proserpio for making me aware of the extensive work on crystal nets that has appeared in the chemistry literature, as well as for valuable comments about the history of nets research and the notation for nets used by crystallographers. I am also indebted to the anonymous referees for a number of helpful suggestions that have improved the article.

 \end{document}